\documentstyle{amsppt}
\magnification1100
\pagewidth{13.5cm}

\define \noi {\noindent}
\define \I{\text {\sl Id}}
\define \Z{\Bbb Z}
\define \R{\Bbb R}
\define \N{\Bbb N}
\define \Q{\Bbb Q}

\define \mpb{\medpagebreak}
\define \spb{\smallpagebreak}

\define \ovr{\overset \to}
\define \undr{\underset \to}
\define \ld{\lambda}
\define \Ld{\Lambda}
\define \g{\gamma}
\define \G{\Gamma}
\define \bt{\beta}
\define \al{\alpha}
\define\k{\kern-15pt} 

\define\phanm{\phantom{$-$}}
\define\phanu{\phantom{$_\frac 12,$}}
\TagsOnRight

\rightheadtext{LENGTH SPECTRA OF FLAT MANIFOLDS}
\leftheadtext{R.J.\ Miatello and J.P.\ Rossetti}
\topmatter
\title  LENGTH SPECTRA AND $P$-SPECTRA OF COMPACT FLAT MANIFOLDS \endtitle
\author R.J.\ Miatello  and J.P.\ Rossetti  \endauthor
\address FaMAF-ciem, Universidad
Nacional de C\'ordoba, 5000 C\'ordoba, Argentina\. 
\endaddress
\email miatello\@mate.uncor.edu \endemail
\address FaMAF-ciem, Univ.\ Nac.\ Cba.\ Argentina.  Currently: Department of Mathematics,  
Dartmouth College, Hanover, NH 03755, USA\. 
\endaddress
\email rossetti\@mate.uncor.edu.  juan.p.rossetti\@dartmouth.edu \endemail

\thanks 2000 {\it Mathematics Subject Classification.} Primary 58J53; \,Secondary 58C22, 20H15. \endthanks
\thanks{Partially supported by Conicet and Secyt U.N.C.}\endthanks
\keywords closed geodesic, length spectrum, isospectral, Bieberbach group \endkeywords

\abstract
We study the length, weak length and complex length spectrum of closed geodesics 
of a compact flat Riemannian manifold, comparing  length-isospectrality with
isospectrality of the Laplacian acting on $p$-forms. 
Using integral roots of the Krawtchouk polynomials, we give many pairs of $p$-isospectral flat manifolds 
having different lengths of closed geodesics and in some cases, different injectivity radius
and different first eigenvalue.

We prove a Poisson summation formula relating the $p$-eigenvalue spectrum with the lengths of closed geodesics. 
As a consequence we show that the spectrum determines the lengths 
of closed geodesics and, by an example, that it does not determine the complex lengths. Furthermore we show that
orientability is an audible property for flat manifolds.
We give a variety of examples, for instance, a pair of
isospectral (resp.\ Sunada isospectral)  manifolds with different length spectra and a pair with the same
complex length spectra and  not $p$-isospectral for any $p$, 
or else  $p$-isospectral for only one value of $p\ne 0$.
\endabstract
\endtopmatter

\document

\heading{Introduction}\endheading

The $p$-spectrum of a compact Riemannian manifold  $M$ is  the collection
of eigenvalues, with multiplicities,  of the Laplacian acting on $p$-forms.
It will be denoted by $\text{spec}_p(M)$.
Two manifolds $M$ and $M'$ are said to be
$p$-isospectral ---or isospectral on $p$-forms--- if $\text{spec}_p(M)=\text{spec}_p(M')$.
The word isospectral is reserved for the function case, i.e.\ it corresponds  to  $0$-isospectral.

Let $\Gamma$ be the fundamental group of $M$. It is well known that the free homotopy classes  of 
closed paths in $M$  are in a one to one correspondence with the conjugacy classes in $\Gamma$.
Furthermore, in each such free homotopy class there is at least a 
closed (i.e.\ periodic) geodesic ---namely the closed path
of smallest length in the class.  In the case when the sectional curvature of $M$ is nonpositive, 
if two closed geodesics are freely homotopic then
they can be deformed into each other by means of a smooth homotopy through a flat 
surface in $M$, hence they  have the same length $l$.  
This length is called the {\it length} of $\gamma$, denoted $l(\g)$, 
where $\gamma \in \Gamma$ is any representative of this class. 
The {\it complex length} of $\g$ is the pair $l_c(\g):=(l(\g),[V])$, where $V \in O(n-1)$ is determined by the 
holonomy of $\g$ (see Section 2) and $[V]$ denotes the conjugacy class.
The {\it multiplicity} of a length $l$ 
(resp.\ of a complex length $(l,[V])$) is defined to be the  number of free homotopy classes having length~$l$ (resp.\
$(l,[V])$).  The  {\it weak length spectrum} (resp.\ {\it weak complex length spectrum}) of $M$, denoted $L$-spectrum
(resp.\
$L_c$-spectrum),  is defined as the set of all lengths (resp.\ complex lengths) of closed geodesics in
$M$, while the {\it length spectrum} (resp.\ {\it complex length spectrum}), denoted  {\it
$[L]$-spectrum} (resp.\ {\it $[L_c]$-spectrum}), is the set of lengths (resp.\ complex lengths) 
of closed geodesics,  with multiplicities. Two manifolds are said to be 
{\it $[L]$-isospectral} or {\it length isospectral} 
if they have the same $[L]$-spectrum. 
(The notions of $L$, $L_c$, and $[L_c]$-isospectrality are defined similarly.) 
They are said to be {\it marked length isospectral} if there
exists a length-preserving isomorphism between their fundamental groups.  

The relationship between length spectrum and eigenvalue spectrum of $M$ has been studied for some time. 
For flat tori and for Riemann surfaces, it is  known that the length spectrum and the eigenvalue spectrum
determine each other (see \cite{Hu1,2}).
Also, it has been proved that ``generically" $\text{spec}(M)$ determines the length spectrum 
of $M$ (see \cite{CdV}). In \cite{DG}, Introduction, an asymptotic formula (see (4.6)) is stated that indicates that an
analogous result holds for $\text{spec}_p(M)$, for any $p\ge 0$. 

All the known
examples of isospectral compact Riemannian manifolds are $L$-isos\-pec\-tral.   In \cite{Go}, C.\ Gordon gave the first
example of a pair of Riemannian manifolds ---they are Heisenberg manifolds--- that are isospectral but not
$[L]$-isospectral. R.\ Gornet in \cite{Gt1,2} gave, among other illuminating examples,  the first example of pairs of
manifolds ---they are 3-step nilmanifolds--- having the same marked length spectrum,   isospectral but not
$1$-isospectral.  For other recent work on the length spectra of nilmanifolds see
\cite{GoM} and \cite{GtM}.  The complex length spectrum has been considered in (\cite{Re}, \cite{Me}) for 
hyperbolic manifolds of dimension $n=3$ and in \cite{Sa} in the case of locally symmetric spaces of negative
curvature.

The goal of this paper is to study the various length spectra for compact flat Riemannian 
manifolds (flat manifolds for short) and to compare the different notions of  length isospectrality with
$p$-isospectrality, for $p\ge 0$. 

We will determine the complex lengths of closed geodesics for general flat manifolds.
We will see that, in general, the $p$-spectrum  does not determine the weak length spectrum. 
For manifolds of diagonal type (see Definition 1.2) this can happen only when  $K^n_p(x)$, 
the {\it (binary) Krawtchouk polynomial of degree} $n$ has integral roots. Using
such roots, we give many pairs of $p$-isospectral flat manifolds  having different lengths of closed
geodesics and in some cases, different injectivity radius (Ex.\ 2.3(i)-(vi)). 
These seem to be the first such examples in the context of compact Riemannian manifolds. 
These examples might be considered quite odd since they do not follow the generic behavior 
and seem to contradict the wave trace formula. 
An explanation on how they are consistent with the
heat  and wave trace formulas is given in Remark 4.11.  

We give several pairs, most of them with different fundamental groups (Examples 3.3 through 3.7), 
comparing length isospectrality with other types of isospectrality.
The examples are obtained by an elementary construction ---they are flat tori  of low dimensions ($n\ge 4$),
divided by free actions of
$\Z_2^r,\, r \le 3$, $\Z_4$ or $\Z_4\times \Z_2$--- and  in particular, it is quite easy to compute their length
spectra and real cohomology. For instance, we give a pair of isospectral, not Sunada isospectral manifolds 
that are not $[L]$-isospectral  (Ex.\ 3.6) and a pair of 4-dimensional flat manifolds that are Sunada
isospectral (see Remark 1.3) but not $[L]$-isospectral (Ex.\ 3.4).  This last example seemed unlikely 
to exist in the context of flat manifolds.

Example 3.5 shows a pair of manifolds both having the same complex length spectrum and not 
isospectral to each other. They can be chosen so that they are not $p$-isospectral for any $p$, or else, 
isospectral for some $p>0$. Examples 3.4 and 3.5 show clear differences with the situation for hyperbolic manifolds, 
where such examples cannot exist (see \cite{GoM}, \cite{Sa}). There are other examples summarized in Table 3.9.
Example 3.8 shows pairs of flat manifolds having the same lengths and/or complex lengths of 
closed geodesics but which are very different from each other; for instance, manifolds having different dimension, 
or so that one of them is orientable and the other not. 
In a different direction, we prove that two flat manifolds having the same marked length spectrum are 
necessarily isometric (Prop.\ 3.10).

One of the ways to connect the eigenvalue spectrum with the length spectrum is via  Poisson summation formulas
or via the Selberg trace formula (see \cite{Bl}, \cite{CdV}, \cite{Hu1,2}, \cite{Pe} and \cite{Su} for
instance). In the case of flat manifolds,  Sunada gave one such formula  in the function case (see \cite{Su}) and
as a consequence he showed that if two flat Riemannian manifolds are isospectral, then the corresponding tori
must be isospectral.  In Section 4 we give a Poisson summation formula for vector bundles 
that is related but is different from Sunada's. In the proof we use the formula 
for multiplicities of eigenvalues obtained in \cite{MR2, Theorem 3.1}. 
As a consequence, we show that the spectrum determines the $L$-spectrum.
Example 3.6 shows that the spectrum does not determine the $L_c$-spectrum.

An open general question is whether orientability is an {\it audible} property, that is, whether isospectral Riemannian 
manifolds should be both orientable or both nonorientable.  Using the Poisson formula we show this question has a
positive answer for flat manifolds. This is not the case for $p$-isospectral flat manifolds; in \cite{MR2} we give 
several pairs of  flat manifolds that are $p$-isospectral for only some values of $p$, one of them orientable and the
other not.  Also, P.\ B\'erard and D.\ Webb (\cite{BW}) have constructed pairs of 0-isospectral surfaces with boundary
that are Neumann isospectral, but not Dirichlet isospectral,  one of them orientable and the other not.
Using the Poisson formula, we show that $p$-isospectrality
of two flat manifolds for some $p\ge 0$ (or $\tau$-isospectrality, for any representation $\tau$ of $O(n)$),
implies that the corresponding tori are isospectral and, furthermore, that the orders of the holonomy groups are 
the same. This is a natural extension of the result of Sunada for 0-isospectral flat manifolds.

Another application of the formula is concerned with flat manifolds of diagonal type. This is a restricted family
but  is still a rich and useful class. For instance, all of the examples constructed in sections 2 and 3, with the
exception of Ex.\ 3.6 and  Ex.\ 2.3(iii), are of diagonal type (see also the examples in \cite{MR1,2,3}). 
Bieberbach groups  in this class are more manageable, for instance it is quite straightforward to compute combinatorially
all the  Betti numbers of the associated manifold. 
We show that  for manifolds of diagonal type, isospectrality ---and also $p$-isospectrality when $K^n_p(x)$ has no 
integral roots--- implies Sunada isospectrality, hence $q$-isospectrality for every $q$.  This extends to all $n$, a
result proved by very different methods in
\cite{MR3} for $n\le 8$.

An outline of the paper is as follows. In Section 1 we recall briefly
some basic facts on Bieberbach groups  and some results from \cite{MR2}. 
In Section 2 we 
give a formula for the complex lengths of closed geodesics and show that, if $p>0$, the $p$-spectrum does not determine the lengths of closed geodesics. 
In Section 3 we
state a criterion for $[L]$ and $[L_c]$-isospectrality,
giving several illustrative examples and counterexamples 
together with a table showing many different possibilities. 
Section 4 is devoted to the Poisson summation formula and the consequences described above.

The authors wish to thank Carolyn Gordon for very useful comments on the contents of this paper.
The second author was at Dartmouth College while part of this paper
was done and would like to thank the great hospitality of the Department of Mathematics,  
specially of Carolyn Gordon and David Webb.

\heading{\S 1 Preliminaries}\endheading

We shall first recall some standard facts on flat Riemannian manifolds (see \cite{Ch}).
A discrete, cocompact subgroup  $\Gamma$ of  the isometry group of $\R^n$,  $I(\R^n)$,
is called a crystallographic group. If furthermore, $\G$ is torsion-free, then $\G$ is said to be a Bieberbach group.
Such $\Gamma$ acts properly discontinuously on $\R^n$, thus $M_\Gamma =
\Gamma\backslash\R^n$ is a compact flat Riemannian manifold with fundamental group $\Gamma$. Any such manifold
arises in this way.  Any element $\gamma \in
I(\R^n)$ decomposes uniquely $\gamma = B L_b$, with $B \in O(n)$ and $b\in \R^n.$ 
The  translations in $\Gamma$
form a normal, maximal abelian subgroup  of finite index, $L_\Lambda$,  $\Lambda$  a lattice in $\R^n$ which is $B$-stable for each $BL_b \in \Gamma$. The quotient 
$F := \Lambda\backslash\Gamma$ is called the holonomy group of $\Gamma$ and gives the linear holonomy group of the Riemannian manifold $M_\Gamma$.
The action of $F$ on $\Lambda$ defines an integral representation of $F$, usually called the  holonomy representation. 
Denote $n_B=\text{dim} \ker(B-\I)$. If $BL_b$ in $\G$, then it is known that $n_B>0$. 
The next lemma contains some facts that we will need. 

\proclaim{Lemma 1.1}
Let $\Gamma$ be a Bieberbach group, and let $\gamma=BL_b \in \Gamma$. Let $L$ be {\it any} lattice 
stable by $B$ and let  $p_B$ denote the orthogonal projection onto $\text{ker}(B-\I)$.
Then we have
\roster
\item"(i)"  $p_B(b) \ne 0$. 
Furthermore $\text{ker} (B-\I)=\text{Im}(B-\I)^\perp$.
\item"(ii)" ${{L}^{\text{B}}}:= L \cap \text{ker} (B-\I)$ is a lattice in $\text{ker} (B-\I)$.
\item"(iii)" Let $L^*= \{\mu \in \R^n : \langle \mu, \lambda \rangle \in \Z, \text{ for any  } \lambda \in L\}$,
the dual lattice of $L$. Then 
${\left({\left({L}^{\text{B}}\right)}^* \right)}^{\text{B}}=p_B({L}^*)$.
\endroster
\endproclaim

\demo{Proof} If $B$ has order $m$, then we have that $(BL_b)^m =L_{Cb}$, where $C=\sum_{j=0}^{m-1} B^j$. 
Since $B\in O(n)$,  then $C=C^t$. Furthermore, since $(B-\I)C=B^m-\I=0$ and $B$ is of order $m$, it follows  that   
$\text{ker} (B-\I) = \text{Im}\, C=({\text{ker}\,C})^{\perp}$. 
Now  $(BL_b)^m \ne \I$, so $Cb\ne 0$, hence  $b\notin 
\text{ker} (B-\I)^{\perp}$, as claimed. 
Now $\text{Im}(B-\I)^\perp=\text{ker} {(B-\I)}^*=\text{ker} (B-\I)$, hence 
the second assertion in (i) is clear.
 
To verify (ii) we note that if $L_\Q=\Q$-span$(L)$, then
$B\, L_\Q= L_\Q$ 
and the $\Q$-rank of the matrix of $B-\I$ on a $\Z$-basis of $L$ equals the $\R$-rank, hence $\text{dim}_\Q \text{ker} (B-\I)\cap {L_\Q}=\text{dim}_\R \text{ker} (B-\I).$ Thus, if $\{v_j:1\le j\le r\}$ is a $\Q$-basis of 
$\text{ker} (B-\I)\cap{L_\Q}$ 
and if  $m_1,\dots,m_r\in \Z\setminus\{0\}$ are such that $m_j v_j \in L,$ for $1\le j\le r$, then $\text{ker} (B-\I) \cap L$ contains $\sum_{j=1}^r \Z m_jv_j$, a lattice in 
$\text{ker} (B-\I)$. This implies the assertion.

Relative to (iii), we set $W=\text{ker} (B-\I)$. 
Since $W^*=W^\perp$ we have:
$$(L \cap W)^* \cap W= ({L}^* + W^\perp) \cap W =  p_W ({L}^*).\qed$$
\enddemo

We now recall from \cite{MR2,3}  some facts on the spectrum of Laplacian operators 
on vector bundles over flat manifolds. If $\tau$ is an irreducible representation of $K=O(n)$ and $G=I(\R^n)$
we form the vector bundle $E_\tau$ over $G/K \simeq \R^n$ associated to $\tau$ and consider the corresponding
bundle $\Gamma\backslash E_\tau$ over $\Gamma\backslash  \R^n =M_\Gamma$. Let $-\Delta_\tau$ be the connection
Laplacian on this bundle.  For any $\mu$ a nonnegative
real number, let $\Lambda^*_\mu=\{\lambda \in \Lambda^* : \| \lambda\|^2=\mu\}$.  In \cite{MR3, Thm.\ 2.1} we
have shown that the multiplicity of the eigenvalue $4\pi^2 \mu$ of 
$-\Delta_\tau$ is given by
$$ d_{\tau,\mu}(\Gamma)= |F|^{-1} \sum_{\gamma=BL_b \in \Lambda\backslash \Gamma}\text{tr}\,\tau(B) \,e_{\mu,\gamma}\tag{1.1}$$
where $e_{\mu,\gamma} = \sum_{v\in {\Lambda^*_\mu}:Bv=v} e^{-2\pi i v.b}$.   
In the case when $\tau=\tau_p$, the $p$-exterior representation of $O(n)$, we shall write $\text{tr}_p(B)$ and
$d_{p,\mu}(\Gamma)$ in place of
$\text{tr}\,\tau_p(B)$ and $ d_{\tau_p,\mu}(\Gamma)$ respectively. 

\mpb

For a special class of flat manifolds the terms in this formula can be made more explicit.

\definition{Definition 1.2} \cite{MR3, Def.\ 1.3.}
We  say that a Bieberbach group $\Gamma$ is of {\it diagonal type} if 
there exists an orthonormal $\Z$-basis $\{e_1,\dots,e_n\}$ of the lattice $\Lambda$ 
such that for any element $BL_b\in\Gamma$, $Be_i=\pm e_i$ for $1\le i\le n$.
Similarly, $M_\G$ is said to be of {\it diagonal type}, if $\G$ is so.
We note that it may be assumed that 
the lattice $\Lambda$ of $\Gamma$ is the canonical lattice.
\enddefinition

These manifolds have, in particular, holonomy group $F\simeq \Z_2^r$, for some $r \le n-1$.
After conjugation by a translation we may assume furthermore that  $b\in \frac 12 \Ld$, for 
any $BL_b \in \G$ (see \cite{MR3, Lemma 1.4}). 
In this case we have the following expressions for the terms $e_{\mu,\gamma}$ in the multiplicity formula (1.1):
$$e_{\mu,\gamma}= \sum_{v\in \Lambda_\mu : Bv=v}\, (-1)^{|I_v^{odd} \cap I_{2b}^{odd}|},\tag1.2$$
where $I_v=\{j:v.e_j \ne 0\}$ and $I_v^{odd}=\{j:v.e_j \text{ is odd}\}$. 
Furthermore, the traces $\text{tr}_p(B)$ are given by integral values of the so called Krawtchouk 
polynomials $K_p^n(x)$ (see \cite{MR2, Remark 3.6} and also \cite{MR3}; see \cite{KL} for more information on Krawtchouk
polynomials).
Indeed, if $n_B= \text{dim} \ker(B-\I)$, we have:
$$\text{tr}_p(B)=K_p^n(n-n_B), \quad \text{ where } K_p^n(x):= \sum_{t=0}^p (-1)^t \binom xt \binom {n-x}{p-t}.\tag{1.3}$$

We shall also use the notation 
$I_B:=\big\{ 1\le i \le n : Be_i=e_i \big\}\,$, so $|I_B|=n_B$  and   
$I_B\cap I_{2b}^{odd}=\{ i\in I_B : b.e_i \equiv \frac 12 \mod \Z \}$. 
We set
$$c_{d,t}(\G) :=\big|\big\{BL_b \in F:n_B=d \text{ and }|I_B\cap I_{2b}^{odd}|=t\big\}\big|,\, \text{ for }
0\le t \le d \le n.\tag{1.4}
$$
We note that by Lemma 1.1 $|I_B\cap I_{2b}^{odd}|>0$, for any $BL_b \in \G$, except when $B=\I$, $b\in \Ld$.

\example{Remark 1.3} In \cite{MR1,3} we gave combinatorial expressions for  
the numbers $c_{d,t}$, called {\it Sunada numbers}. We showed that their equality  
for $\G$ and $\G'$ is equivalent to have that $M_\G$ and $M_{\G'}$ verify the conditions in Sunada's theorem, 
that is, they are {\it Sunada isospectral} (see \cite{MR3} Def.\ 3.2, Thm.\ 3.3 and the discussion following it).  
In particular $c_{d,t}(\G)=c_{d,t}(\G')$  for every $d,t$ implies 
that $M_\G$ and $M_{\G'}$ are $p$-isospectral for all $p$.
A different proof of this fact will be given in \S 4 (Theorem 4.5).
\endexample

\heading{\S 2 Length of closed geodesics and $p$-spectra}\endheading

Let $M$ be a compact Riemannian manifold  of nonpositive curvature. 
If $\al (t)$ is a closed geodesic on $M$ with period $t_o$, the parallel transport $\tau$ along
$\al(t)$ from $\al(0)$ to $\al(t_o)=\al(0)$ is such that $\tau (\dot{\al}(0))=\dot \al(0)$, hence it defines an
element
$V
\in O((\R \dot\al(0))^\perp)\simeq O(n-1)$. One can associate to $V$ a well defined conjugacy class $[V]$ in
$O(n-1)$,  called the {\it holonomy} of $\al(t)$.   By definition, the complex length of $\al$ is the pair
$l_c(\al):=(l(\al),[V])$; $l_c(\al)$ depends only on the free homotopy class of $\al$ and not on $\al$.

Now given $l\ge 0$ 
let $m(l)$  denote the multiplicity of the length $l$, 
that is, the number of free homotopy classes of closed geodesics 
in $M$ such that $l(\al)=l$. 
It is known that these multiplicities are finite.

The $L$, $L_c$, $[L]$ and {\it $[L_c]$-spectrum}
of $M$ and the respective notions of isospectrality have been defined in the Introduction.

Clearly if $M$ and $M'$ are $L_c$-isospectral (resp.\ $[L_c]$-isospectral) then they 
are $L$-isospectral (resp.\ $[L]$-isospectral).

We will start by proving a proposition that gives some basic properties 
of closed geodesics in a flat manifold $M_\Gamma$.
Let $\Gamma$ be an $n$-dimensional Bieberbach group and $\g =BL_b\in \G$.
For any $v \in \R^n$ we may write 
$$v=v_+ + v' \,\,\text{ with }\,\, v_+ \in
\text{ker}(B-\I) \,\,\text{ and }\,\, v' \in
\text{ker}(B-\I)^\perp.\tag2.1$$ We note that Lemma 1.1(i) says that $b_+\ne 0$. 
We take (2.1) as the definition of $v_+$ and~$v'$.

We have that $\|BL_{b}x - x\|^2= \|(B-\I)x+ Bb' + b_+\|^2=\|(B-\I)x+ Bb'\|^2 +\| b_+\|^2,$ by Lemma 1.1(i). 
Since $B-\I$ is an isomorphism when restricted to $\text{ker}(B-\I)^\perp$, then $x$ can always be chosen so that the
first summand is zero.  This says that $\inf\big\{\text{dist}(\gamma x,x): x\in \R^n\big\}$ can be attained at some $z$
and it equals $\| b_+\|$. Let $o_\g$ be defined uniquely by: 
$$(B-\I)o_\g = - Bb', \quad o_\g\in \text{ker} (B-\I)^\perp.\tag2.2$$

Now, if $z$ is such that $\|BL_{b}z - z\|=\|b_+\|$, then the first summand in the above expression for 
$\|BL_{b}z -z\|^2$ is zero, so we see that $z$ is of the form $o_\g+ u$ for some $u \in \text{ker}(B-\I)$.

Actually,  the next proposition allows to characterize  the points $z \in \R^n$ such that
$\text{dist}(\gamma z,z)=\inf\big\{\text{dist}(\gamma x,x): x\in \R^n\big\}= \| b_+\|$ as the points lying on lines on $\R^n$ stable by $\g$. These points form
the affine space $o_\g+\ker(B-I)$ of dimension $n_B$.

For $u\in\ker(B-\I)$, $t\in\R$, 
we define  
$\,\,\al_{\g,u}(t):= o_\g + u + t b_+$, . We set 
$\,\al_{\g}(t):= o_\g +t b_+$, $t\in\R$, i.e. $\al_\g=\al_{\g,0}$.

\proclaim{Proposition 2.1} Let $\Gamma$ be a Bieberbach group. 
\roster
\item"(i)" If $\g =BL_b\in \G$, then $\g$ preserves the lines $o_\g + u + \R b_+$, $o_\g$ and $u$ as above.  
Furthermore $\g \al_{\g,u}(t)=\al_{\g,u}(t+1)$.
Any line in $\R^n$ stable by $\g$ is of the form $o_\g + u + \R b_+$, for some $u \in \text{ker} (B-\I)$.  
\item"(ii)" The geodesic $\al_{\g,u}$ pushes down to a closed geodesic $\bar\alpha_{\g,u}(t)$ in $M_\G$, 
$t\in\Z\backslash\R$, of length $l(\bar \al_{\g,u})=\| b_+\|$.  
Any closed geodesic in $M_\G$ is of the form $\bar \al_{\g,u}$ for some $\g=BL_b \in \G$ and $u \in \text{ker}(B-\I)$. 
\item"(iii)" $\bar \al_{\g,u}$  is freely homotopic to $\bar \al_\g$. 
The holonomy of $\bar \al_{\g,u}$ is given by $[B^\perp]$, where $B^\perp$ denote the restriction of $B$ to 
$(\R b_+)^\perp$. 
\item"(iv)" The $L$-spectrum (resp.\ $L_c$-spectrum) of $M_\Gamma$ is the set of numbers $\|b_+\|$ (resp.\ the set of pairs $(\|b_+\|, [B^\perp])$), where  $BL_b$ runs through all elements of $\Gamma$.  
\item"(v)" The $[L]$-spectrum (resp.\ $[L_c]$-spectrum) of $M_\Gamma$ is the set of numbers
$\|b_+\|$ (resp.\ the set of pairs $(\|b_+\|, [B^\perp])$), counted with multiplicities, where $\gamma = BL_{b}$ runs through a full set of representatives for the $\Gamma$-conjugacy classes in $\Gamma$. 
\endroster
\endproclaim
\demo{Proof}
If $\g=BL_b$ and $t\in \R$,  taking into account the definition of $o_\g$, we have
$$\g \al_{\g,u}(t)=b_+ + Bb' +B o_\g +u + tb_+ = o_\g+ u + (t+1)b_+=\al_{\g,u}(t+1).$$

Now let $w+\R v$ be a line stable by $\g$. This happens if and only if
\,\,$Bb+ Bw +\R Bv= w+\R v$, or equivalently
$$\R Bv = \R v \quad\text{ and }\quad b_+ +Bb'+(B-\I)w \in \R v.$$ 
It follows that $B v=\pm v$, since $B\in O(n)$. 
Since $b_+ \ne 0$ and $Bb'+(B-\I)w \in \text{ker} (B-\I)^\perp$ we have that $v_+\ne 0$, thus $v=v_+$.  
Hence  $Bv=v$ and furthermore $\R v = \R b_+$ and $Bb' + (B-\I)w =0$. 
Now, by the definition of $o_\g \in \text{ker} (B-\I)^\perp$, 
then clearly  \,$u := w-o_\g\in \text{ker} (B-\I)$. 
This implies the second assertion in (i).

\spb 

For (ii), it is clear by (i) that $\bar \al_{\g,u}(t)$ 
is a closed geodesic in $M_\G$, with length equal to the length 
of the segment in $\R^n$ from $\al_{\g,u}(0)$ to $\al_{\g,u}(1)$, 
which equals~$\|b_+\|$.
Since any closed geodesic in $M_\G$ is the push down of a geodesic in $\R^n$ that is translated into itself by
some $\g \in \G$, then the second assertion in (ii) 
follows from (i).
\spb

Relative to (iii), we see that $\bar \al_{\g,su}$ with $s\in [0,1]$ is a continuous family of closed geodesics 
in $M_\G$, which shows that $\bar \al_{\g,u}$ and $\bar \al_{\g}=\bar \al_{\g,0}$ are freely homotopic. 

To determine $l_c(\g)$, we note that the parallel transport 
along $\bar\al_{\g,u}$, from $\pi(o_\g +u)$ to itself, 
is given by $B$ and since $Bb_+ =b_+$, then $B$ preserves $(\R b_+)^\perp$. 
This implies that the holonomy of $\g$ is $[B^\perp]$. 

Finally, assertions (iv) and (v) follow immediately from (ii),(iii). 
\qed 
\enddemo

In the sequel 
we shall denote by  diag$(a_1,\dots, a_n)$, the $n\times n$
diagonal matrix with $a_j$ in the $j^{\text{th}}$ diagonal entry.

\example{Example 2.2} {\it The Klein bottle.}  
As a warm-up, we will first look at the simplest case of the Klein-bottle group. 
We will determine the conjugacy classes in
$\G$ and  the closed geodesics, computing their lengths and their
respective multiplicities. 
We let $\G=\langle BL_b,\Ld\rangle$ with $B=\text{diag}(-1,1)$, 
$b= \frac {e_2}2$, $\Ld = \Z e_1 + \Z e_2$. Thus, $M_\G$ is a flat Klein bottle. We have 
$$\G= \{L_\ld : \ld \in \Ld \} \cup \{BL_{b+\ld} : \ld \in \Ld \},$$
a disjoint union.
We first compute the conjugacy classes in $\G$. We have:
$$BL_b L_\ld {(BL_{b})}^{-1} = L_{B\ld},\qquad 
L_\mu BL_b L_{-\mu}= BL_{b+(B-Id)\mu},$$
for any $\ld, \mu \in \Ld$. Thus 
$$L_{m_1e_1+m_2e_2} \sim L_{-m_1e_1+m_2e_2},\qquad
BL_b \sim BL_{b + 2ke_1},$$
for any $m_1,m_2,k \in \Z$, where $\sim$ means $\G$-conjugate.

Thus, a full set of representatives for the $\G$-conjugacy classes is 
$$ \{L_{m_1e_1+m_2e_2} : m_1 \in \N_0, m_2 \in \Z\} \cup  \{BL_b L_{m_1e_1+m_2e_2} : m_1=0,1, m_2 \in \Z\}.$$
The corresponding lengths are given by:
$$ l(L_{m_1e_1+m_2e_2})=
{\left( m_1^2 +m_2^2 \right)}^{\frac{1}2},\qquad\quad l(BL_{b+m_1e+m_2e_2})=|\tfrac{1}2
+m_2|.$$
Now, if $x \in \R^2, \g =BL_{b+\ld} \in \G$, with $\ld =m_1e_1+m_2e_2$, 
then the segment in $\R^2$ joining $x$ to $\g x= Bx +B\ld +\frac{e_2}2$ has length
$\| \g x -x\|={\left(\|(B-\I)x - m_1e_1\|^2+(\frac 12+m_2)^2\right)}^\frac 12 $.
This segment pushes down to a closed (periodic) geodesic in $M_\G$ if and only if it has 
minimal length, that is, $x=- \frac {m_1}2 e_1+ se_2$, $s\in \R$, and the length is 
$l(\g)=\big|m_2+\frac 12\big|$. Thus,
the closed geodesics in $M_\G$ are the pushdowns of the vertical segments joining 
$-\frac {m_1}2 e_1+se_2$ and $\frac {-m_1}2 e_1+(\frac 12+m_2 +s) e_2$, for any $m_1,m_2 \in \Z$, $s\in [0,1)$, 
together with the pushdowns of the segments joining $x$ to $x+\ld$, for any $x \in \R^2$ and any $\ld \in \Ld$. 

We may also see this in a different way, by noticing that   
by Proposition 2.1, $BL_{b+\ld}$ stabilizes the line
$\al_\g(t)=-\frac {m_1}2 e_1 + t(m_2+ \frac 12)e_2$, since $o_{BL_{b+\ld}}= -\frac {m_1}2 e_1$.
Hence $\bar\al_\g(t)$, the pushdown of $\al_\g(t)$, is a closed geodesic in $M_\G$. 
We note that $L_{e_1}\al_{\g}(t) = \al_{\g{L_{-2e_1}}}(t)$, for all $t$, hence 
$\bar \al_{\g} = \bar\al_{\g{L_{-2e_1}}}$, thus we may assume that $m_1 \in \{0,1\}$.

Using the above set of representatives for the $\G$-conjugacy classes we
see that a set of representatives for the free homotopy classes of closed paths are the pushdowns of the segments 
joining $(0,0)$ to $\ld \in \Ld$ with $m_1\ge 0$,
together with the pushdowns of the segments joining  
$(0,0)$ to $ (m_2+ \frac 12){e_2}$ and those joining 
$\frac {e_1}2$ to $\frac {e_1}2 +(\frac 12+m_2) {e_2}$, $m_2\in\Z$.

The multiplicities of the lengths of closed geodesics are as follows. 
If $l^2\in\N_0$ then $m(l)$ equals the number of solutions of the equation $m_1^2+ m_2^2 = l^2$, with $m_1\in \N_0$, $m_2
\in \Z$. 
On the other hand, if $l =\frac 12 +k$, $k\in \N_0$, we find that $m(l)=4$. Indeed, the conjugacy
classes with length $\frac 12 +k$ correspond to the segments joining, either $(0,0)$ 
with $\pm(\frac 12+k)e_2$, or $\frac{e_1}2$ with $\frac{e_1}2\pm(\frac 12+k)e_2$. 

Regarding the holonomy of $\g$, it equals 1 for $\g$ such that 
${l(\g)}^2 \in \N_0$, i.e.\ $\g \in \Ld$, and it equals 
$-1$ if  $l(\g)\in \tfrac 12 +\N_0$.  

We note that a more general Klein-bottle group $\G_{\al,\bt}$ for, $\al,\bt \in \R_{>0}$, can be defined by taking $\Ld :=\Ld_{\al,\bt}= \Z\al e_1 + \Z \bt e_2$, $b=\frac \bt 2  e_2$ and with $B$ as before (see \cite{BGM}). The lengths of closed geodesics in $M_{\G_{\al,\bt}}$ are in this case  given by $|m_2+\frac 12|\bt$ and $({m_1}^2\al^2 + {m_2}^2 \bt^2)^{\frac 12}$, for 
$m_1, m_2 \in \Z$.
\endexample

\example{Example 2.3} The goal of this example is to show that the $p$-spectrum for some $0<p\le n$ 
does not determine the lengths of closed geodesics.
We will construct
several pairs of $p$-isospectral manifolds having different $L$-spectrum. 
In some cases, the smallest lengths of closed geodesics are distinct for 
$M_\G$ and $M_{\G'}$. Hence, the injectivity radii ($\frac 12$ of the smallest length of closed geodesics) are 
distinct for $M_\G$ and $M_{\G'}$ (see Ex.3(ii),(iii),(iv),(vi)). 
Also, some examples have different first eigenvalue on functions and, among them,  
some have different injectivity radius. 
In Remark 4.11 we explain how these examples are consistent with the heat and wave trace formulas.

\mpb
 
(i) We  take $\G$ and $\G'$ having holonomy group $\Z_2$. 
The nontrivial elements in $F$, $F'$ are given by $BL_b$ with $b=\tfrac {e_1}2$, $B= \text{diag}(1,-1,-1,-1)$
for $\G$ and $B'= \text{diag}(1,1,1,-1)$ for $\G'$.  Then $\G$ and $\G'$ are nonisomorphic Bieberbach groups (see
\cite{MR2, Ex.\ 4.1}). By using formula (1.1), we get that the multiplicities of the eigenvalue $4\pi
^2$ equal
$d_{0,1}(\G)= \tfrac 12 (8 -2)=3$ and 
$d_{0,1}(\G')= \tfrac 12 (8 + (-2+2+2))=5$, hence both manifolds are not isospectral. Furthermore 
$\text{tr}_p( \text{diag}(1,-1,-1,-1))=K_p^4(3)$, $\text{tr}_p( \text{diag}(1,1,1,-1))=K_p^4(1)$ and  one has that 
$K_p^4(3)= K_p^4(1)=0$ if and only if $p=2$ (see (1.3)). 
It follows from (1.1) that $M_\G$ and $M_{\G'}$ are $p$-isospectral 
if and only if $p=2$ (\cite{MR2, Ex.\ 4.1}).

On the other hand, the lengths $l$ such that $l^2 \not \in \Z$ have the form $\big|\frac 12 +m\big|$ for $M_{\G}$ and 
$\left({(\frac 12 +m_1)^2 +m_2^2 + m_3^2}\right)^\frac{1}2$ for $M_{\G'}$, for arbitrary $m, m_i \in \Z$.
This implies that $M_\G$ and $M_{\G'}$ are not $L$-isospectral.
For instance, $\frac {\sqrt{5}}2$ is a length of a closed geodesic for $M_{\G'}$ but not for $M_{\G}$. 

Regarding complex lengths we note for instance that if we take
$\g=BL_b, \g'=B'L_b$ we have $l(\g)=l(\g')=\tfrac 12$, the  minimal length, however the complex lengths 
of $\g,\g'$ are different, since $l_c(\g)=(\tfrac 12, \text{diag}(-1,-1,-1))$ 
and $l_c(\g')=(\tfrac 12, \text{diag}(1,1,-1))$.

\mpb

\noi(ii) We consider a variation of the previous example (see \cite{MR2, Ex.\ 4.2}).
Again $\G$ and $\G'$ both have holonomy groups $\Z_2$ and we let 
$B= \text{diag}(1,1,-1,-1)$ for both $\G$ and  $\G'$, taking
$b=\tfrac {e_1}2$ for $\G$ and $b'=\tfrac {e_1 + e_2}2$ for $\G'$. 
In this case $\G$ and $\G'$ are isomorphic but the corresponding manifolds are not isometric ($\G'$ is obtained by
conjugating
$\G$ by  a $C \in GL(4,\R)$). 
Since  $K_p^4(2)=0$ if and only if
$p=1,3$, it follows that the associated flat manifolds are isospectral only for these values of $p$. 
They are not isospectral nor $L$-isospectral. Indeed, the 
lengths that do not correspond to lattice elements have the form
$\left((\frac 12 +m_1)^2 +m_2^2\right)^\frac{1}2$ 
for $M_\G$, and $\left((\frac 12 +m_1)^2 +(\tfrac 12 + m_2)^2\right)^\frac{1}2$  for $M_{\G'}$,  
respectively, for arbitrary $m_1,m_2 \in \Z$. In particular $\frac 12$ is a
length for $\G$ but not for
$\G'$. Thus $M_\G$ and $M_{\G'}$ have injectivity radius equal to $\frac 14$ and $\frac {\sqrt{2}}4$ respectively.  

\spb

We note that (i) and (ii) can be generalized to any even dimension $n \ge 4$, as done in \cite{MR2}, 
with the same properties. 
In particular, in (ii) we obtain manifolds that are $p$-isospectral for any $p$ odd, $0<p<n$,
and are not $L$-isospectral.

\mpb

(iii)
We now let $\G$ and $\G'$ be as in Example 5.8 in \cite{MR2}. Then $\Gamma$ and $\Gamma'$ are nonisomorphic 
Bieberbach groups and $M_\G$ and $M_{\G'}$ have dimension 4 and holonomy groups $\Z_2^2$ and $\Z_4$ respectively.  
Let $\Gamma = \langle B_1L_{b_1}, B_2L_{b_2}, \Lambda \rangle$, and 
$\Gamma' = \langle B'L_{b'}, \Lambda \rangle$, where 
$\Lambda = \Z^4$, $b_1=\frac{e_1}2,b_2=\frac{e_4}2, b'=\frac{e_4}4$, $B_1=\text{diag}(1,1,-1,-1)$, 
$B_2=\text{diag}(1,-1,-1,1)$ and $B'=\text{diag}({\widetilde J},-1,1)$,  
where $\undr{{\widetilde J}}=\left[\smallmatrix 0&1 \\ -1&0 \endsmallmatrix\right]$. 
 
The following table lists the nontrivial elements in $F$ and $F'$ in a convenient notation, writing the rotational parts in columns, together with subindices that indicate the 
nonzero translational components: for instance, since $b_1=\frac{e_1}2$, we write $\frac{1}2$ as a subindex of
the first diagonal element of $B_1$.

\

\centerline{
\vbox{\tabskip=0pt
\def\tablerule{\noalign{\hrule}}
\halign to142pt {\strut#& \vrule#\tabskip=1em plus 2em&
 \hfil# \hfil&  #& 
 \hfil# \hfil&  #& 
 \hfil# \hfil& \vrule# 
 \tabskip=0pt\cr\tablerule
&& {\vrule height14pt width0pt depth8pt} $B_1$ && $B_2$ && $B_1B_2$  & \cr\tablerule 
&& {\vrule height12pt width0pt depth0pt}  $1_{\frac 12}$  && \phanm$1$\phanu  && \phanm$1_{\frac 12}$  &\cr
&&  \phanm$1$\phanu  && $-1$\phanu  && $-1$\phanu &\cr
&&  $-1$\phanu  && $-1$\phanu  && \phanm$1$\phanu &\cr
&&  $-1\undr{}$\phanu  && \phanm$1_{\frac 12}$  && $-1_{\frac 12}$ &\cr
\tablerule
\cr}}
\qquad
\vbox{\tabskip=0pt
\def\tablerule{\noalign{\hrule}}
\halign to140pt {\strut#& \vrule#\tabskip=1em plus 2em&
 \hfil# \hfil&  #& 
 \hfil# \hfil&  #& 
 \hfil# \hfil& \vrule# 
 \tabskip=0pt\cr\tablerule
&& {\vrule height14pt width0pt depth8pt} $B'$ && ${B'}^2$ && ${B'}^3$  & \cr\tablerule 
&& {\vrule height12pt width0pt depth0pt}  $\ovr{\undr{{\widetilde J}}}$  && $-\I_2$  && $-{\widetilde J}$\phanu  &\cr
&&  $\ovr{\undr{-1}}$\phanu  && \phanm$1$\phanu  && $-1$\phanu &\cr
{\vrule height0pt width0pt depth8pt}&& \phanm$\ovr{1_{\frac 14}}$ && \phanm$1_{\frac 1{\undr 2}}$ && \phanm$1_{\frac 34}$  &\cr
\tablerule
\cr}}
}

We have that $\text{tr}_p (\I)= \binom 4p$ and, for the three nontrivial 
elements $B$ in $F$, $\text{tr}_p (B)=K_p^4(2)$ 
equals $0$ if $p=1,3$, it equals $1$ if $p=0,4$ and it equals $-2$ if $p=2$. 

In the case of $\G'$ there are some differences. We have 
$n_{B'}=n_{{B'}^3}=1$ and $n_{{B'}^2}=2$. Also, $b'=b'_+=\tfrac{e_4}4$.
For $C:=B', {B'}^3$, we have that $\text{tr}_p (C)=0$ if $p=1,2,3$, and 
$\text{tr}_p (C)=1$, if $p=0,4$. Furthermore, $\text{tr}_p ({B'}^2)=K_p^4(2)$. Using this information we get
that $M_\G$ and
$M_{\G'}$ are $p$-isospectral if and only if $p=1,3$. Also,  $M_\G$ is orientable, while $M_{\G'}$ is not.
They are not $L$-isospectral, since $M_{\G'}$ has closed geodesics with length $\frac 14$ , while $M_{\G}$ does not. 
Actually, they have different injectivity radius.

\mpb

(iv) Let $\G$ and $\G'$ of dimension 4, with holonomy groups $\Z_2^2$.
Set $B_1=B'_1= \text{diag}(1,1,-1,-1)$, $b_1=\tfrac {e_1+e_2+e_3+e_4}2$, $b'_1=\tfrac {e_1+e_2}2$;  
$B_2=B'_2= \text{diag}(1,1,-1,1)$, $b_2=\tfrac {e_4}2$, $b'_2=\tfrac {e_2}2$. 
The lattice is $\Z^4$.
By a verification of the conditions (i) and (ii) in \cite{MR2, Prop.\ 2.1},
one has that $\Gamma$ and $\Gamma'$ are nonisomorphic Bieberbach groups.

It is not difficult to see, by an argument similar to that in (i)-(iii), that 
these manifolds are 2-isospectral, they are not $L$-isospectral, and the first (nonzero) 
eigenvalue is $4\pi^2$ for $M_{\G'}$ and $8\pi^2$, for $M_\G$.

\spb

These manifolds have the same injectivity radius; but if we change
in $\G$ the value of $b_2$ into $\tfrac {e_2+e_4}2$, then we obtain 
manifolds having the same spectral properties as before but now they do not have 
the same injectivity radius, namely $\frac {\sqrt{2}}4$ and $\frac 14$ respectively. 
In column notation the last groups are as follows.

\

\centerline{
\vbox{\tabskip=0pt
\def\tablerule{\noalign{\hrule}}
\halign to142pt {\strut#& \vrule#\tabskip=1em plus 2em&
 \hfil# \hfil&  #& 
 \hfil# \hfil&  #& 
 \hfil# \hfil& \vrule# 
 \tabskip=0pt\cr\tablerule
&& {\vrule height14pt width0pt depth8pt} $B_1$ && $B_2$ && $B_3$  & \cr\tablerule 
&& {\vrule height12pt width0pt depth0pt}  \,$1_{\frac 12}$  && \phanm$1$\phanu  && \phanm$1_{\frac 12}$  &\cr
&&  \phanm$1_{\frac 12}$  && \phanm$1_{\frac 12}$  && \phanm$1$\phanu &\cr
&&  $-1_{\frac 12}$  && $-1$\phanu  && \phanm$1_{\frac 12}$ &\cr
&&  $-1_{\frac 12}$ && \phanm$\undr1_{\frac 12}$ &&  $-1$\phanu &\cr
\tablerule
\cr}}
\qquad
\vbox{\tabskip=0pt
\def\tablerule{\noalign{\hrule}}
\halign to150pt {\strut#& \vrule#\tabskip=1em plus 2em&
 \hfil# \hfil&  #& 
 \hfil# \hfil&  #& 
 \hfil# \hfil& \vrule# 
 \tabskip=0pt\cr\tablerule
&& {\vrule height14pt width0pt depth8pt} $B'_1$ && $B'_2$ && $B'_3$  & \cr\tablerule 
&& {\vrule height12pt width0pt depth0pt}  \,$1_{\frac 12}$  && \phanm$1$\phanu  && \phanm$1_{\frac 12}$  &\cr
&&  \phanm$1_{\frac 12}$   && \phanm$1_{\frac 12}$  && \phanm$1$\phanu &\cr
&&  $-1$\phanu && $-1$\phanu  &&  \phanm$1$\phanu &\cr
&&  $-1$\phanu && \phanm$1$\phanu  && $-1$\phanu  &\cr
\tablerule
\cr}}
}

\mpb

(v) For each value of $p$ one can give pairs $\G$, $\G'$ of dimension $n=2p$, with holonomy group $\Z_2$, 
that are $q$-isospectral only for $q=p$  and such that they are not $L$-isospectral.
More generally, following Proposition 3.15 in \cite{MR3}, it is possible to construct, for each integral zero of
the Krawtchouk polynomial, namely $K^n_p(j)=0$ with $j\ne \tfrac n2$, groups $\G$, $\G'$ of dimension $n$, with holonomy group $\Z_2$, that are $p$ and $(n-p)$-isospectral and not $L$-isospectral.
They are $q$-isospectral if and only if $K_q^n(j)=0$ and this happens generically,
if and only if $q=p$ or $q=n-p$.

\mpb

(vi)
Each example 
in (v) can be extended to give large families of $p$-isospectral 
flat manifolds of dimension $n$, having pairwise different $L$-spectrum. We sketch this construction in a particular case,
since all cases are similar. Take for $1\le j\le k<n$, $n$ even, the Bie\-ber\-bach groups  
$\G^n_{k,j}:=\langle  C_kL_{\frac{e_1+\dots +e_j}2}, \Z^n \rangle$, where 
$C_k:= \text{diag} (\underbrace{1,1,\dots,1}_k,\underbrace{-1,\dots,-1}_{n-k})$. 
For fixed $k$, these groups are isomorphic but the corresponding manifolds are not
isospectral. For $k$ odd, they are $\frac n2$-isospectral (see \cite{MR2, Ex.\ 4.2} for more details).
One has that 
if $j_1<j_2$, both $\not\equiv 0\mod 4$, then $\G^n_{k,j_1}$ and $\G^n_{k,j_2}$ are not L-isospectral, since 
$\frac {\sqrt{j_1}}2$ is a length for $\G^n_{k,j_1}$ but not for $\G^n_{k,j_2}$.
For fixed $k$, the family $\G^n_{k,j}$ with $1\le j\le k$, $j\not\equiv 0 \mod 4$ has cardinality 
approximately equal to $\frac {3k}4$. The corresponding manifolds are all $p$-isospectral to each other for $p=\frac n2$, but they have different lengths of closed geodesics. 
Observe that it is possible to take $k=n-1$, giving a family of cardinality approximately equal to $\frac{3n}4$.
\endexample

\heading{\S 3 Length spectrum and  $p$-spectra}\endheading

In this section we will consider the question of $[L]$-isospectrality of
flat manifolds. The next criterion will be useful.
If $\bt\in \G , \G'$ let $[\bt]$ denote the conjugacy class of $\bt$.

\proclaim{Proposition 3.1} Let $\G, \G'$ be Bieberbach groups with translation lattices $\Ld, \Ld'$ respectively. 
Suppose there exist partitions $\Cal P$ and
$\Cal P'$ of $\Ld\backslash\G$ and $\Ld'\backslash\G'$ respectively and a 
bijection $\phi : \Cal P \rightarrow \Cal P'$
such that for every $c\in\R_{>0}$, $O\in O(n-1)$ and $P \in \Cal P$, the cardinality of
$\,\,\left\{ [\g L_\ld]\colon \ld \in\Ld,\,\g \in P,\, 
l(\g L_\ld)=c \,\, (\text{resp. } l_c(\g L_\ld)=(c,[O]))\right\}$ 
equals the cardinality of
$\,\,\left\{ [\g'L_{\ld'}]\colon\ld'\in\Ld',\g'\in \phi(P),\, 
l(\g' L_{\ld'})=c \,\,(\text{resp. }
l_c(\g'L_{\ld'})=(c,[O])) \right\},$
then $\G\backslash\R^n$ and $\G'\backslash\R^n$ are length (resp.\ complex length) isospectral.
\endproclaim

When applying this criterion in Examples 3.3 and 3.7 below we shall use the point partition, that is, 
each class in $\Cal P$ and $\Cal P'$ will have one
element, hence $\phi$ will be a bijection from $\Ld\backslash\G$ to $\Ld'\backslash\G'$.  Example 3.5 has less standard
spectral properties and will require a less obvious partition of $F$ and~$F'$.

In order to be able to compute the $[L]$-spectrum and $[L_c]$-spectrum 
of a general flat manifold $M_\G$ one needs a parametrization of the conjugacy classes of $\G$. 
This is in general complicated
but it becomes much simpler when the Bieberbach group $\G$  is of diagonal type.
For a general $\G$, we have that if $\g_i=B_iL_{b_i}, \g_j=B_jL_{b_j} \in \G, \ld, \mu  \in \Ld$, conjugation of $L_\ld$ by $\g_i$, and $\g_i L_\ld $ by $L_\mu$, yield the following relations:
$$\align
 L_\ld &\sim L_{B_i \ld},\qquad
 \g_i L_\ld \sim \g_i L_\ld L_{(B_i^{-1}-Id)\mu}.
\tag{3.1}
\endalign
$$
Furthermore, if the holonomy group is abelian, conjugation of $\g_iL_\ld$ by $\g_j$, $j\ne i$, yields
$$
\align
B_j L_{b_j} B_i L_{b_i+\ld} L_{-b_j} B_j^{-1}& =  B_i L_{B_j( (B_i^{-1}-Id)b_j +b_i+\ld)}\\
&= B_i L_{b_i+\ld}L_{B_j(B_i^{-1}-Id)b_j + (B_j-Id)(b_i+\ld)}\\
&=B_i L_{b_i+B_j\ld}L_{B_j(B_i^{-1}-Id)b_j + (B_j-Id)b_i}.
\endalign
$$
We thus get:
$$
 \g_i L_\ld \sim \g_i L_{B_j\ld} L_{(B_j-Id)b_i + B_j(B_i^{-1}-Id)b_j},\text{ for } j \ne i.
\tag{3.2}
$$

\example{Remark 3.2}
We now mention some simple facts that are rather direct consequences of (3.1), (3.2) and the definitions.

\mpb

(i) We note that $\lambda_1 , \lambda_2 \in \Lambda$ are $\Gamma$-conjugate if and only if 
there is $BL_b \in \Gamma$ such that $B\lambda_1 = \lambda_2$.  If $\G$ and $\G'$ have the same 
lattices and the same integral holonomy representations, then   the multiplicities of the lengths corresponding to
lattice elements are the same. (In many of the examples these are exactly the lengths $l$ such that $l^2\in\N$.)

\mpb

(ii) If $M_\G$ and $M_{\G'}$ are Sunada isospectral then they are $L_c$-isospectral.

\spb

Indeed, if $\beta=CL_c \in I(\R^n)$ and $\g=BL_b \in \G$ then
$\bt \g \bt^{-1}= CBC^{-1}L_{C((B^{-1}-Id)c+b)}$
and furthermore,
${(C((B^{-1}-\I)c+b))}_+=Cb_+$. 
Thus, $l(\bt \g \bt^{-1})=\|Cb_+\|=\|b_+\|$. The holonomy component of $\bt \g \bt^{-1}$ is $\left[
{CBC^{-1}}_{|{(\R Cb_+)^\perp}}\right]=[B_{|{(\R b_+)}^\perp}]$. 
Hence, this  shows that $l_c(\bt \g \bt^{-1})=l_c(\g)$. Note that this also shows that the complex length 
is an invariant of the conjugacy class of $\g$ in $I(\R^n)$, in particular the complex length spectrum is
well defined.

Thus, using that $l_c(\bt \g L_\ld \bt^{-1})=l_c(\g L_\ld)$
it follows that $M_\G$ and $M_\G'$ are  $L_c$-isospectral.
 
\mpb

(iii) In the notation of Prop.\ 2.1, assume that for each $P \in {\Cal P}$ (resp.\ ${\Cal P}'$) and 
for any $\g_1, \g_2 \in P$ 
the holonomy components of $\g_1, \g_2$ are the same. 
Suppose that there exists $\phi$ as in Prop.\ 2.1 satisfying the conditions for $[L]$-isospectrality and suppose also 
that $\phi$ preserves holonomy components. 
Then $\G\backslash\R^n$ and $\G'\backslash\R^n$ are $[L_c]$-isospectral.
\endexample

\example{Example 3.3} We let $\Gamma$ and $\Gamma'$ as in  \cite{MR2, Ex. 4.5}. In column notation:

\

\centerline{
\vbox{\tabskip=0pt
\def\tablerule{\noalign{\hrule}}
\halign to146pt {\strut#& \vrule#\tabskip=1em plus 2em&
 \hfil# \hfil&  #& 
 \hfil# \hfil&  #& 
 \hfil# \hfil& \vrule# 
 \tabskip=0pt\cr\tablerule
&& {\vrule height14pt width0pt depth8pt} $B_1$ && $B_2$ && $B_1B_2$  & \cr\tablerule 
&& {\vrule height12pt width0pt depth0pt}  \,$1$\phanu  && \phanm$1$\phanu  && \phanm$1$\phanu  &\cr
&&  \phanm$1_{\frac 12}$  && $-1$\phanu  && $-1_{\frac 12}$ &\cr
&&  $-1$\phanu  && \phanm$1_{\frac 12}$  && $-1_{\frac 12}$ &\cr
&&  $\undr -1_{\frac 12}$  && $-1$\phanu  && \phanm$1_{\frac 12}$ &\cr
\tablerule
\cr}}
\qquad
\vbox{\tabskip=0pt
\def\tablerule{\noalign{\hrule}}
\halign to146pt {\strut#& \vrule#\tabskip=1em plus 2em&
 \hfil# \hfil&  #& 
 \hfil# \hfil&  #& 
 \hfil# \hfil& \vrule# 
 \tabskip=0pt\cr\tablerule
&& {\vrule height14pt width0pt depth8pt} $B'_1$ && $B'_2$ && $B'_1 B'_2$  & \cr\tablerule 
&& {\vrule height12pt width0pt depth0pt}  \,$1$\phanu  && \phanm$1_{\frac 12}$  && \phanm$1_{\frac 12}$  &\cr
&&  \phanm$1_{\frac 12}$  && $-1$\phanu  && $-1_{\frac 12}$ &\cr
&&  $-1$\phanu  && \phanm$1$\phanu  && $-1$\phanu &\cr
&&  $\undr -1$\phanu  && $-1$\phanu  && \phanm$1$\phanu &\cr
\tablerule
\cr}}
}

Then $\Gamma$ and $\Gamma'$ both have the same holonomy representation of diagonal type, with holonomy group $\Z^2_2$. 
In {\cite{MR2} it was shown that $M_\G$ and $M_{\G'}$ are isospectral,
actually they are Sunada isospectral but they are not diffeomorphic. They are $[L]$-isospectral and actually $[L_c]$-isospectral, as we  shall see by using Proposition 3.1. 

We take $\Cal P$ the point partition and the bijection $\phi$
as the identity. Since $B_i=B'_i$ for each $i$, $\phi$ preserves the holonomy components of the complex lengths $l_c(\g)$, $\g \in \G$.

Since the lattices and the integral holonomy representations are the same for both manifolds then by Remark 3.2(i), 
the multiplicities of the lengths of the elements of the form $L_\ld$ are the same for $\G$ and $\G'$.

The remaining lengths are of the form $l=\left({{(\frac 12+k_1)}^2+k_2^2}\right)^\frac{1}2$ where $k_1,k_2\in\Z$.
In $\G$ the elements of length $l$ with rotational part $B_1$ are of the form $B_1 L_{b_1+\ld}$ with $\ld=m_1e_1+m_2e_2+m_3e_3+m_4e_4$, 
$m_i\in\Z$ for $1\le i\le 4$ and
such that $m_1^2+(\frac 12+m_2)^2=l^2$.
The second relation in (3.1) implies that 
$B_1 L_{b_1+\ld} \sim B_1 L_{b_1+\ld+2\delta_3 e_3+2\delta_4 e_4}$ for $\delta_3,\delta_4\in\Z$,
thus we may assume that $m_3,m_4\in\{0,1\}$ without leaving out any conjugacy class. 
Indeed (3.1) implies that $m_3$ and $m_4$ can be taken modulo $2$ in this case.

Now by (3.2) $B_1 L_{b_1+\ld} \sim B_1 L_{b_1+m_1e_1-(m_2+1)e_2+(m_3-1)e_3-(m_4+1)e_4}$, and there are no other relations among elements with rotational part  $B_1$.
Since the second coordinate has the form $\frac 12 +m_2$ 
and changes into $-\frac 12 -m_2$ when applying (3.2), this implies that in each conjugacy class in $\G$ there are exactly two elements related in this way. That is, among the elements 
$B_1 L_{v}$ in $\G$ with $v=m_1e_1+(\frac 12+m_2)e_2+m_3e_3+(\frac 12+m_4)e_4$, 
$m_1,m_2\in\Z$, $m_3,m_4=0$ or $1$, and such that $m_1^2+(\frac 12+m_2)^2=l^2$, the number of conjugacy classes 
is exactly half the total number of elements.

If we proceed to do the same calculation in $\G'$ for the elements with rotational part $B'_1=B_1$,
we see that the set of elements with length $l$ are exactly the same as in the previous case; also (3.1) gives 
the same relations as before, and now by (3.2) the relations are 
$B_1 L_{b'_1+\ld} \sim B_1 L_{b'_1+m_1e_1-(m_2+1)e_2+m_3e_3+m_4e_4}$.
As before we see that the elements of a fixed length $l$ occurring here are divided by a
factor or 2 when we apply (3.2) to take conjugacy classes.  This proves the equality of cardinalities as required in
Proposition 3.1 for the pair $B_1$ and $\phi(B_1)$.

It is not difficult to check that 
the elements with rotational parts $B_2$ and $B_3$ 
can be handled in a completely similar way. After applying (3.1), by (3.2), the elements 
go in pairs to form a conjugacy class 
and in all cases the number of classes with a given length 
corresponding to $B_i$, $i=2$ or $i=3$,  is exactly the same for both $\G$ and $\G'$.

Now the above discussion implies that the conditions in Remark 3.2(iii) are satisfied, hence both manifolds are $[L_c]$-isospectral.

\endexample

\example{Example 3.4} We now consider a simple pair of Bieberbach groups $\G$, $\G'$, of dimension 4, 
with holonomy group $\Z_2^2$ and of diagonal type. 
We shall see that $M_\G$ and $M_{\G'}$ are Sunada isospectral, but not $[L]$-isospectral. Let
$\G=\langle B_1L_{b_1}, B_2L_{b_2}, \Ld \rangle$,
$\G'=\langle B_1L_{b'_1}, B_2L_{b'_2}, \Ld \rangle$,
let $B_3=B_1B_2$, where $B_i,b_i,b'_i$ are given in the following table.

\

\centerline{
\vbox{\tabskip=0pt
\def\tablerule{\noalign{\hrule}}
\halign to142pt {\strut#& \vrule#\tabskip=1em plus 2em&
 \hfil# \hfil&  #& 
 \hfil# \hfil&  #& 
 \hfil# \hfil& \vrule# 
 \tabskip=0pt\cr\tablerule
&& {\vrule height14pt width0pt depth8pt} $B_1$ && $B_2$ && $B_3$  & \cr\tablerule 
&& {\vrule height12pt width0pt depth0pt}  \,$1_{\frac 12}$  && \phanm$\undr 1_{\frac 12}$  && \phanm$1$\phanu  &\cr
&&  \phanm$1$\phanu  && \phanm$1_{\frac 12}$  && \phanm$\undr1_{\frac 12}$ &\cr
&&  \phanm$1$\phanu  && $-1$\phanu  && $-1$\phanu &\cr
&&  $-1$\phanu  && \phanm$1$\phanu  && $-1$\phanu &\cr
\tablerule
\cr}}
\qquad
\vbox{\tabskip=0pt
\def\tablerule{\noalign{\hrule}}
\halign to150pt {\strut#& \vrule#\tabskip=1em plus 2em&
 \hfil# \hfil&  #& 
 \hfil# \hfil&  #& 
 \hfil# \hfil& \vrule# 
 \tabskip=0pt\cr\tablerule
&& {\vrule height14pt width0pt depth8pt} $B'_1$ && $B'_2$ && $B'_3$  & \cr\tablerule 
&& {\vrule height12pt width0pt depth0pt}  \,$\undr1$\phanu  && \phanm$1$\phanu  && \phanm$1$\phanu  &\cr
&&  \phanm$1$\phanu  && \phanm$1_{\frac 12}$  && \phanm$1_{\frac 12}$ &\cr
&&  \phanm$1_{\frac 12}$  && $-1_{\frac 12}$  && $-1$\phanu &\cr
&&  $-1_{\frac 12}$  && \phanm$\undr1_{\frac 12}$  && $-1$\phanu &\cr
\tablerule
\cr}}
}

The corresponding manifolds are  Sunada isospectral (see Theorem 4.5 and Remark
1.3) since the numbers $c_{d,t}$ , $0\le t \le d\le 4$, are the same for $\G, \G'$. Indeed, 
$c_{2,1}=c_{3,1}=c_{3,2}=c_{4,0}=1$ and $c_{d,t}=0$ for the other values of $d,t$. 

The lengths of closed geodesics corresponding to nonlattice elements 
are of the form:
${\left((\frac 12 + m_1)^2 +m_2^2+m_3^2 \right)}^\frac{1}2\,$ and
${\left( (\frac 12 +m_1)^2+(\frac 12 +m_2)^2+m_3^2\right)}^\frac{1}2\,$ with 
$m_i \in \Z, 1\le i \le 3$, for both $\G$ and $\G'$.
Furthermore, since the manifolds are Sunada isospectral, then the $L_c$-spectra are the same by Remark 3.2(ii).
We will now see that $M_\G$ and $M_{\G'}$ are not
$[L]$-isospectral, by showing that the multiplicity of the length $\frac 12$ is different for $M_\G$ and $M_{\G'}$.

Among the elements in $\G$ having length $\frac 12$ some are of the form $B_1L_{\frac {e_1}2+\ld}$ with $\ld =-m_1e_1+m_4e_4$, where $m_1 \in \{0,1\}, m_4 \in \Z$.
Since, by (3.1), $B_1L_{\frac {e_1}2} \sim B_1L_{\frac {e_1}2+ 2ke_4}$ for any $k \in \Z$, then we may also take $m_4\in \{0,1\}$.
We easily see that relations (3.1), (3.2) do not give any identifications between these elements, 
so we get 4 different conjugacy classes with length $\frac 12$.

The other elements of length $\frac 12$ are 
$B_3L_{\frac {e_2}2+\ld}$ with $\ld =-m_2e_2+m_3e_3+m_4e_4$, where $m_i \in \{0,1\}$ (here we again use relation (3.1) to have $m_3,m_4 \in \{0,1\}$).
Furthermore, by (3.2) we have
$$B_3 L_{\frac {e_2}2+\ld} \sim B_3L_{\frac {e_2}2+B_j\ld}=B_3L_{\frac {e_2}2 -m_2e_2\pm m_3e_3\pm m_4e_4}$$
hence (3.2) gives no new relations and  these elements lie in 8 different conjugacy classes. Thus the length $l= \frac 12$ has multiplicity 12 in $M_\G$.

For $\G'$ the elements of length $\frac 12$ have the form
$B'_1L_{\frac {e_3+e_4}2+\ld}$ with $\ld =-m_3e_3+m_4e_4$, $m_3,m_4 \in \{0,1\}$ and $B'_3L_{\frac {e_2}2+\ld}$ with $\ld =-m_2e_2+m_3e_3+m_4e_4$, $m_i\in \{0,1\}$. We have by (3.2) the following relations:
$$
\align
B'_1L_{\frac {e_3+e_4}2+\ld}\sim& B'_1L_{\frac {e_3+e_4}2+B'_2 \ld}
L_{(B'_2-Id)\frac{e_3+e_4}2 +B'_2(B'_1-Id)\frac{e_2+e_3+e_4}2}\\
=&B'_1L_{\frac {e_3+e_4}2+(m_3-1)e_3 +(m_4-1)e_4}
\endalign
$$
This implies that these 4 elements determine 2 conjugacy classes in $\G'$.  
Similarly one computes that the remaining 8 elements for $B'_3$ determine 4 conjugacy classes in $\G'$. 
Hence the length $l=\frac 12$ has multiplicity 6 in $M_{\G'}$.

Examples of manifolds with similar spectral properties are given in \cite{Go, Ex. 2.4(a)} and in \cite{Gt1, Ex.I}, 
by using 2 and 3-step nilmanifolds, respectively.
We note that such an example cannot exist for hyperbolic manifolds since strongly isospectral 
implies $[L]$-isospectral in this context (see \cite{GoM}).
\endexample

\example{Example 3.5} 
We will now see that, in the context of flat manifolds, the
$[L]$-spectrum (and even the $[L_c]$-spectrum) does not determine the $p$-spectrum for any $p$, $0\le p\le n$.  
This shows a difference with the case of hyperbolic manifolds, 
since in this context $[L_c]$-isospectral implies strongly isospectral (see \cite{Sa}).
Indeed, we
will construct two flat manifolds of dimension 13 (resp.\ 14), with holonomy  group isomorphic to $\Z_2^3$, which are
$[L_c]$-isospectral but are not $p$-isospectral for any $0\le p\le 13$ (resp.\ for any $0\le p\le 14$, except for $p=7$). 
The corresponding Bieberbach groups have both the canonical lattice and the same integral holonomy representation.  
The translational parts of elements in $\G$ with nontrivial rotational part differ from the corresponding elements
in $\G'$ only in the last nine coordinates, where the rotational parts act as the identity.  

We shall represent the elements of $F$, $F'$ in column notation, writing on the right (resp.\ left) the coordinates corresponding to elements in $\G$ (resp.\ $\G'$).

\

\centerline{
\vbox{\tabskip=0pt
\def\tablerule{\noalign{\hrule}}
\halign to374pt {\strut#& \vrule#\tabskip=1em plus 2em&
 \hfil# \hfil&  #& 
 \hfil# \hfil&  #& 
 \hfil# \hfil&  \vrule#& 
 \hfil# \hfil&  #& 
 \hfil# \hfil&  #&  
 \hfil# \hfil&  #& 
 \hfil# \hfil&  \vrule#  
\tabskip=0pt\cr\tablerule
&& {\vrule height14pt width0pt depth8pt} $B_1$ && $B_2$ && $B_3$ && $B_1B_2$  && $B_1B_3$ && $B_2B_3$ &&  $B_1B_2B_3$ & \cr\tablerule 
&& {\vrule height12pt width0pt depth0pt}  $1$  && $-1$\phanu  && $-1$\phanu && $-1$\phanu && $-1$\phanu && $1$ && $1$ &\cr
&& $-1$\phanu  && $1$  && $-1$\phanu && $-1$\phanu && $1$ && $-1$\phanu && $1$ &\cr
{\vrule height0pt width0pt depth8pt}&& $-_{\frac 12}1_{\frac 12}$\phanu && $-1$\phanu && $1$ && $_{\frac 12}1_{\frac 12}$ && $-_{\frac 12}1_{\frac 12}$\phanu && $-1$\phanu &&  $_{\frac 12}1_{\frac 12}$ &\cr
{\vrule height0pt width0pt depth8pt}&& $-1$\phanu && $-1$\phanu && $-1$\phanu && $1$ && $1$ && $1$ && $-1$\phanu &\cr
{\vrule height0pt width0pt depth8pt}&& $_{\frac 12}1_{\frac 12}$ && $_{\frac 12}1$\phanu && $1$ && \phanm$1_{\frac 12}$ && $_{\frac 12}1_{\frac 12}$ && $_{\frac 12}1$\phanu && \phanu$1_{\frac 12}$ &\cr
{\vrule height0pt width0pt depth8pt}&& $_{\frac 12}1$\phanu && $_{\frac 12}1_{\frac 12}$ && $1$ && \phanu$1_{\frac 12}$ && $_{\frac 12}1$\phanu && $_{\frac 12}1_{\frac 12}$ && \phanu$1_{\frac 12}$ &\cr
{\vrule height0pt width0pt depth8pt}&& $_{\frac 12}1$\phanu && $_{\frac 12}1_{\frac 12}$ && \phanu$1_{\frac 12}$ && \phanu$1_{\frac 12}$ && $_{\frac 12}1_{\frac 12}$ && $_{\frac 12}1$\phanu && $1$ &\cr
{\vrule height0pt width0pt depth8pt}&& $_{\frac 12}1_{\frac 12}$ && $_{\frac 12}1$\phanu && $_{\frac 12}1_{\frac 12}$ && \phanu$1_{\frac 12}$ && $1$ && \phanu$1_{\frac 12}$ && $_{\frac 12}1$\phanu &\cr
{\vrule height0pt width0pt depth8pt}&& $_{\frac 12}1$\phanu && $1$ && \phanu$1_{\frac 12}$ && $_{\frac 12}1$\phanu && $_{\frac 12}1_{\frac 12}$ && \phanu$1_{\frac 12}$ && $_{\frac 12}1_{\frac 12}$ &\cr
{\vrule height0pt width0pt depth8pt}&& $_{\frac 12}1$\phanu && $1$ && \phanu$1_{\frac 12}$ && $_{\frac 12}1$\phanu && $_{\frac 12}1_{\frac 12}$ && \phanu$1_{\frac 12}$ && $_{\frac 12}1_{\frac 12}$ &\cr
{\vrule height0pt width0pt depth8pt}&& $1$ && $_{\frac 12}1$\phanu && \phanu$1_{\frac 12}$ && $_{\frac 12}1$\phanu && \phanu$1_{\frac 12}$ && $_{\frac 12}1_{\frac 12}$ && $_{\frac 12}1_{\frac 12}$ &\cr
{\vrule height0pt width0pt depth8pt}&& $1$ && $_{\frac 12}1$\phanu && \phanu$1_{\frac 12}$ && $_{\frac 12}1$\phanu && \phanu$1_{\frac 12}$ && $_{\frac 12}1_{\frac 12}$ && $_{\frac 12}1_{\frac 12}$ &\cr
{\vrule height0pt width0pt depth8pt}&& $1$ && $1$ && $_{\frac 12}1$\phanu && $1$ && $_{\frac 12}1$\phanu && 
$_{\frac 12}1$\phanu && $_{\frac 12}1$ &\cr
\tablerule
\cr}}
}

\

One sees that $\G$ and $\G'$ are not isospectral. 
Indeed, one verifies that the contribution of $\g_1=B_1L_{b_1}$ to the multiplicity formula (1.1)   
is different from the contribution of $\g'_1=B_1L_{b'_1}$  
for $\G'$.  
On the other hand, the total contribution of the remaining elements is the same for both $\G$ and $\G'$, 
as we can see by taking the bijection $\phi:B_2L_{b_2} \leftrightarrow B_3L_{b'_3}$,  $B_3L_{b_3} 
\leftrightarrow B_2L_{b'_2}$ and $\phi=\I$ for the remaining 
elements.

We now check the $[L]$-isospectrality by applying the criterion in Proposition 3.1. 
We shall use the partitions ${\Cal P}$, ${\Cal P}'$ 
of $\Ld\backslash \G$ and $\Ld\backslash \G'$ respectively, such that any class in ${\Cal P}$ and ${\Cal P'}$ has exactly one element except for
$\{\g_1,\g_2\} \in {\Cal P}$ and $\{\g_1',\g_2'\} \in {\Cal P}'$.
We choose the bijection $\phi$ so that it maps $\{\g_1,\g_2\}$ to $\{\g_3'\}$,
$\{\g_3\}$ to $\{\g_1',\g_2'\}$ and it equals the identity, otherwise. 
We must show that for this choice of $\phi$ the conditions in Proposition 3.1 are satisfied.

It is clear that the lengths 
corresponding to elements of the form $L_\ld$ have the same multiplicities in 
both manifolds by Remark 3.2(i).
The same is true for the elements of the form 
$BL_{b+\ld}$ with $B=B_1B_2, B_1B_3, B_2B_3$ and $B_1B_2B_3$, since they play 
exactly the same role in $\G$ and in $\G'$. 

We claim that the multiplicities of the lengths, when restricted to elements of 
the form $B_3L_{b_3 +\ld}$, are the same as the combined multiplicities of the 
lengths when restricted to elements of the form $B_1L_{b'_1 +\ld}$ and  $B_2L_{b'_2 +\ld}$, $\ld\in\Ld$. 
A similar statement can be made when we compare the combined contributions to the length spectrum of 
elements in $\G$ of the form $B_1L_{b_1 +\ld}$ and  $B_2L_{b_2 +\ld}$ with that of the elements $B_3L_{b'_3 + \ld}$ in $\G'$.

Recall now that the elements of the form $B_3L_{b_3+\ld}$ with 
$\Ld\ni\ld=(m_1,m_2,\dots,m_{13})$ have length 
$$
{\left(m_3^2+m_5^2+m_6^2+{(\tfrac 12+m_7)}^2+{(\tfrac 
12+m_{8})}^2+\dots+{(\tfrac 12+m_{12})}^2+ m_{13}^2\right)}^\frac 
12.\tag3.3$$ 
By conjugating by $L_{e_i}$ for $i=1,2,4$, we may assume that 
$m_1,m_2,m_4\in \{0,1\}$. 
According to (3.2) the only extra relation among these elements 
is given by $B_3L_{b_3+\ld}\sim B_3L_{b_3+\nu}$ where $\nu$ equals $\ld$ except for  
a sign change in only one coordinate, namely 
$\nu=(m_1,m_2,-m_3,m_4,m_5,m_6,\dots,m_{13})$. 
Meanwhile, in 
$\Gamma'$, $B_1L_{b'_1+\ld}$ and $B_2L_{b'_2+\ld}$, as $\ld\in\Ld$ varies, have 
the same lengths as in (3.3). Here we are using that the holonomy group is abelian.
For $B_1L_{b'_1+\ld}$ we may assume $m_2,m_3,m_4\in \{0,1\}$. The remaining 
relations among these elements are:  $B_1L_{b'_1+\ld}\sim B_1L_{b'_1+\nu}$ 
where $\nu$ equals $\ld$ except for a sign change in the first  
coordinate;  $B_1L_{b'_1+\ld}\sim B_1L_{b'_1+\nu'}$ where 
$\nu'=(-m_1,m_2,1-m_3,m_4,m_5,\dots,m_{13})$; and 
the composition of these two relations which gives  $B_1L_{b'_1+\ld}\sim B_1L_{b'_1+\nu''}$ where 
$\nu'' =(m_1,m_2,1-m_3,m_4,m_5,\dots,m_{13})$. 

By taking into account all these relations, one can check that the 
multiplicities of a length for the elements of the form $B_3L_{b_3+ \ld}$ equal 
twice the multiplicities of the same length for the elements of the form $B_1 L_{b'_1+\ld}$.
The same is true for 
$B_2L_{b'_2+\ld}$ in place of $B_1L_{b'_1+\ld}$. Hence, the contribution of 
$B_3L_{b_3+\ld}$ to the $[L]$-spectrum of $M_\G$, as $\ld\in\Ld$ varies, turns out 
to be the same as the contribution of $B_1L_{b'_1+\ld}$ and 
$B_2L_{b'_2+\ld}$ to the $[L]$-spectrum of $M_{\G'}$, as $\ld\in\Ld$ varies.

Similarly, the same happens with $B_3L_{b'_3+\ld}$ when compared with 
$B_1L_{b'_1+\ld}$ and $B_2L_{b'_2+\ld}$ taken together. Therefore $M_\G$ and 
$M_{\G'}$ are $[L]$-isospectral. Furthermore
the conditions in Remark 3.2(iii) are satisfied, thus
$M_\G$ and $M_{\G'}$ are actually $[L_c]$-isospectral.

It is easy to check that $M_\G$ and $M_{\G'}$ are not $p$-isospectral for any $p$, by using
Theorem 3.6(ii), since the coefficients $K^{13}_p(3)\ne 0$ for every $p$.

Finally, if we modify a little bit these manifolds enlarging them to dimension 14, we can obtain
an example of a pair of manifolds with the same properties as before with the only exception that they become $7$-isospectral.  
The change is achieved by replacing the fourth row by two rows, one of the form $(-1,-1,1,1,-1,-1,1)$ 
and the other of the form $(1,1,-1,1,-1,-1,-1)$. Thus $\text{tr}_p(B_j)=K_p^{14}(4)$, for $1\le j \le 3$ 
which is zero if and only if $p=7$. We do not include the verification in this case for brevity.
\endexample

\example{Example 3.6} 
We now briefly consider $\G$ and $\G'$ as in  \cite{MR1, Ex.\ 4.1}.
Here $n=6$ and $F\simeq F' \simeq \Z_4\times \Z_2$. In the standard column notation we may represent the nontrivial
elements in $F, F'$ as follows:

\

\centerline{
\vbox{\tabskip=0pt
\def\tablerule{\noalign{\hrule}}
\halign to340pt {\strut#& \vrule#\tabskip=1em plus 2em&
 \hfil# \hfil&  #& 
 \hfil# \hfil&  #& 
 \hfil# \hfil& \vrule#& 
 \hfil# \hfil& \vrule#& 
 \hfil# \hfil&  #&  
 \hfil# \hfil&  #& 
 \hfil# \hfil&  \vrule#  
\tabskip=0pt\cr\tablerule
&& {\vrule height14pt width0pt depth8pt} $B_1$ && $B_1^2$ && $B_1^3$ && $B_2$  && $B_1B_2$ && $B_1^2B_2$ &&  $B_1^3B_2$ & \cr\tablerule 
&& {\vrule height12pt width0pt depth0pt}  $\widetilde J$  && $-\I_2$  && $-\widetilde J$ && \!\!$-\I_2$ && 
\!\!$-\widetilde J$ && $\I_2$ && $\widetilde J$ &\cr
&&  $\widetilde J$  && $-\I_2$  && $-\widetilde J$ && $\I_2$ && $\widetilde J$ && \!\!$-\I_2$ && \!\!$-\widetilde J$ &\cr
{\vrule height0pt width0pt depth8pt}&& $1_{\frac 14}$\! && $1_{\frac 12}$\! && $1_{\frac 34}$\! && 1 && 
$1_{\frac 14}$\!\! && $1_{\frac 12}$\!\! && $1_{\frac 34}$\!\! &\cr
{\vrule height0pt width0pt depth8pt}&& 1 && 1 && 1 && $1_{\frac 12}$\!\! && $1_{\frac 12}$\!\! && $1_{\frac 12}$\!\! && $1_{\frac 12}$\!\! &\cr
\tablerule
\cr}}
}

\

\centerline{
\vbox{\tabskip=0pt
\def\tablerule{\noalign{\hrule}}
\halign to348pt {\strut#& \vrule#\tabskip=1em plus 2em&
 \hfil# \hfil&  #& 
 \hfil# \hfil&  #& 
 \hfil# \hfil& \vrule#& 
 \hfil# \hfil& \vrule#& 
 \hfil# \hfil&  #&  
 \hfil# \hfil&  #& 
 \hfil# \hfil&  \vrule#  
\tabskip=0pt\cr\tablerule
&& {\vrule height14pt width0pt depth8pt} $B'_1$ && ${B'_1}^2$ && ${B'_1}^3$ && $B'_2$  && $B'_1B'_2$ && ${B'_1}^2B'_2$ &&  ${B'_1}^3B'_2$ & \cr\tablerule 
&& {\vrule height12pt width0pt depth0pt}  $\widetilde J$  && \!$-\I_2$  && \!$-\widetilde J$ && \!$-\I_2$ && 
\!$-\widetilde J$ && $\I_2$ && $\widetilde J$ &\cr
{\vrule height0pt width0pt depth8pt} && 1 && 1  && 1 && \!\!\!$-1$ && \!\!\!$-1$ && \!\!\!$-1$ && \!\!\!$-1$ &\cr
{\vrule height0pt width0pt depth8pt}&& \!\!$-1$ && 1  && \!\!$-1$ && $1_{\frac 12}$\! && \!$-1_{\frac 12}$ && $1_{\frac 12}$\! && \!$-1_{\frac 12}$ &\cr
{\vrule height0pt width0pt depth8pt}&& \!\!$-1$ && $1$ && \!\!$-1$ && 
\!$-1_{\frac 12} $ && $1_{\frac 12}$\! && \!$-1_{\frac 12}$ && 
$1_{\frac 12}$\! &\cr
{\vrule height0pt width0pt depth8pt}&& $1_{\frac 14}$\! && $1_{\frac 12}$\! && $1_{\frac 34}$\! && $1$ && $1_{\frac 14}$\! && $1_{\frac 12}$\! && $1_{\frac 34}$\! &\cr
\tablerule
\cr}}}
\

We have proved in \cite{MR1} that $M_\G$ and $M_{\G'}$ are 0 and 6-isospectral, but not $p$-isospectral for $p\ne 0,6$. We now show they are not $[L]$-isospectral by comparing the multiplicities of the length
$l=\frac 1{\sqrt 2}$. 

The only elements in $\G$ with this length are of the form $B_1^2B_2L_{\frac{e_5+e_6}2 +\ld}$ 
where $\ld=m_3e_3+m_4e_4-m_5e_5 -m_6e_6$, with $m_3,m_4 \in \Z, m_5,m_6 \in \{0,1\}$. 
Using relation (3.1) we may assume that $m_3,m_4\in \{0,1\}$, thus having 16 elements with length $\frac 1{\sqrt{2}}$. 
Now relation (3.2) gives  
$B_1^2B_2L_{\frac{e_5+e_6}2}L_\ld \sim B_1^2B_2L_{\frac{e_5+e_6}2}L_{B\ld}$, with $B\in F$. This implies 
---by taking $B=B_1$--- that the choice 
$m_3=0$, $m_4=1$ is equivalent to the choice $m_3=1$, $m_4=0$, and there are no other identifications.
This gives 12 classes with length $l=\frac 1{\sqrt 2}$ in $\G$.
An entirely similar argument shows that in $\G'$ there are only 8 such classes. This proves the assertion. 

We note also  that $M_\G$ and $M_{\G'}$ have the same $L$-spectrum but they have different $L_c$-spectrum, since for instance, the holonomy component of $\g_1=B_1L_{\frac{e_5}4}$ 
is not the holonomy component of any element in $\G'$.
\endexample

\example{Example 3.7} 
We will describe two 7-dimensional flat manifolds 
of diagonal type which are Sunada isospectral, with isomorphic fundamental groups
(hence diffeomorphic, by Bieberbach's second theorem), 
$[L_c]$-isospectral but not marked length isospectral.

We take
$\G=\langle B_1L_{b_1},B_2L_{b_2},L_\ld \colon \ld\in\Z^7 \rangle$,
$\G'=\langle B'_1L_{b'_1},B'_2L_{b'_2},L_\ld \colon \ld\in\Z^7 \rangle$ as follows

\

\centerline{
\vbox{\tabskip=0pt
\def\tablerule{\noalign{\hrule}}
\halign to146pt {\strut#& \vrule#\tabskip=1em plus 2em&
 \hfil# \hfil&  #& 
 \hfil# \hfil&  #& 
 \hfil# \hfil& \vrule# 
 \tabskip=0pt\cr\tablerule
&& {\vrule height14pt width0pt depth8pt} $B_1$ && $B_2$ && $B_1B_2$  & \cr\tablerule 
&& \phanm$\ovr1_{\frac 12}$  && \phanm$1_{\frac 12}$  && \phanm$1$\phanu  &\cr
&& \phanm$1_{\frac 12}$  && \phanm$1_{\frac 12}$  && \phanm$1$\phanu  &\cr
&&  \phanm$1_{\frac 12}$  && $-1$\phanu  && $-1_{\frac 12}$ &\cr
&&  \phanm$1$\phanu  && $-1$\phanu  && $-1$\phanu &\cr
&&  $-1$\phanu  && \phanm$1_{\frac 12}$  && $-1_{\frac 12}$ &\cr
&&  $-1$\phanu  && \phanm$1$\phanu  && $-1$\phanu &\cr
&&  $-1_{\frac 12}$  &&  $-1$\phanu  && \phanm$\undr1_{\frac 12}$ &\cr
\tablerule
\cr}}
\qquad
\vbox{\tabskip=0pt
\def\tablerule{\noalign{\hrule}}
\halign to152pt {\strut#& \vrule#\tabskip=1em plus 2em&
 \hfil# \hfil&  #& 
 \hfil# \hfil&  #& 
 \hfil# \hfil& \vrule# 
 \tabskip=0pt\cr\tablerule
&& {\vrule height14pt width0pt depth8pt} $B'_1$ && $B'_2$ && $B'_1 B'_2$  & \cr\tablerule 
&&  \phanm$\ovr1_{\frac 12}$  && \phanm$1_{\frac 12}$  && \phanm$1$\phanu  &\cr
&& \phanm$1$\phanu  && \phanm$1$\phanu  && \phanm$1$\phanu  &\cr
&&  \phanm$1_{\frac 12}$  && $-1$\phanu  && $-1_{\frac 12}$ &\cr
&&  \phanm$1_{\frac 12}$  && $-1$\phanu  && $-1_{\frac 12}$ &\cr
&&  $-1$\phanu  && \phanm$1_{\frac 12}$  && $-1_{\frac 12}$ &\cr
&&  $-1$\phanu  && \phanm$1_{\frac 12}$  && $-1_{\frac 12}$ &\cr
&&  $-1_{\frac 12}$  &&  $-1$\phanu  && \phanm$\undr1_{\frac 12}$ &\cr
\tablerule
\cr}}
}

\

By computing the Sunada numbers it is straightforward to check that the manifolds $M_\G$ and $M_{\G'}$ are 
Sunada isospectral (see Remark 1.3).  
Furthermore it is not hard to give an explicit isomorphism from $\G$ to $\G'$ by conjugation by
an affine motion.
  
By using Proposition 3.1 one can show that they are 
indeed $[L_c]$-isospectral, with arguments similar to those in Example 3.3.

We now show that the manifolds are not marked length isospectral and hence not isometric.
Indeed, suppose that 
there exists a length-preserving isomorphism $\phi : \G \to \G'$. Then we must have  $\phi(\Ld)=\Ld$ and also
$\phi(\underset\Z\to{\text{span}}\{ e_1,e_2\})=\underset\Z\to{\text{span}}\{e_1,e_2\}$, 
since this is the space fixed by the action of the holonomy
group.  Furthermore $\phi^{-1}(e_2)=\pm e_i$ with $i=1$ or $2$.
Since $l(B_1L_{b_1})=\frac{\sqrt{3}}2$ and $\phi$ is length-preserving, it follows that $\phi(B_1L_{b_1})=B'_jL_{c'_j}$ where $B'_j=B'_1$ or $B'_2$ and $c'_j=\pm\frac{e_1}2 +w$ 
for some $w\in \frac 12 \underset\Z\to{\text{span}} \{e_3,e_4,\dots,e_7\}$ such that  $\|{c'_j}_+\|=\frac{\sqrt{3}}2$.

Now we have $\,l(B_1L_{b_1+e_i})=\|(b_1+e_i)_+\|=\sqrt{{(\frac 32)}^2+{(\frac 12)}^2+{(\frac 12)}^2}=\frac{\sqrt{11}}2$; 
while 
$\,l(\phi(B_1L_{b_1+e_i}))=l(\phi(B_1L_{b_1})\circ\phi(L_{e_i}))=
l(B'_jL_{b'_j}L_{e_2})=\sqrt{{(\frac 12)}^2+1^2+{(\frac 12)}^2+{(\frac 12)}^2}=\frac{\sqrt{7}}2$,
a contradiction. Hence  $M_\G$ and $M_{\G'}$ are not  marked length isospectral.
\endexample

\example{Example 3.8}
There exist flat manifolds having the same lengths of closed geodesics but which are very different from each other. 
We will now give several $L$-isospectral pairs having, either different dimension, or nonisomorphic
fundamental  groups, or one of them orientable and the other not.
We will make use of the following classical theorem, proved by Lagrange:
{\it Every nonnegative integer can be written as a sum of four squares} (see \cite{Gr} for instance).
As a consequence of this fact we see that {\it all canonical tori $\Z^n\backslash \R^n$, $n\ge 4$, have the same $L$-spectrum}.

The same is true for many other flat manifolds. For instance, if we take 
the Bieberbach groups $\G^n_k:=\langle C_kL_{\frac{e_1}2}, \Z^n\rangle$,  where  
$C_k= \text{diag} (\underbrace{1,1,\dots,1}_k,\underbrace{-1,\dots,-1}_{n-k})$, as in Ex.\ 2.3(vi) 
then, for $5\le k<n$, all groups are $L$-isospectral. This is the case since the 
$(\frac 12+m_1)^2+m_2^2+\dots+m_5^2$, with $m_i\in\Z$, represent 
every number of the form $\frac 14+m$, $m\in\N$.
However, it is clear that $\G^n_k$ are not $L_c$ isospectral nor $[L]$-isospectral to each other.  
 For fixed $k$ and varying $n$, they have different dimensions $n$; half of them are orientable and half
are nonorientable.
\endexample  

\

\topcaption
{\bf Table 3.9.}
\endcaption

\centerline{
\vbox{\tabskip=0pt
\def\tablerule{\noalign{\hrule}}
\halign to378pt {\strut#& \vrule #\tabskip=1.0em plus 0em&
 \hfil# \hfil&  #& 
 \hfil# \hfil&  #& 
 \hfil# \hfil&  #& 
 \hfil# \hfil&  #& 
 \hfil# \hfil&  #&  
 \hfil# \hfil&  \kern-6pt \vrule#  
\tabskip=0pt\cr\tablerule
&& \k $\smallmatrix \text{Some}\\ p\text{-isospectrality?} \endsmallmatrix$ \k
&& \k $\smallmatrix \text{Sunada}\\ \text{isospectral?} \endsmallmatrix$  \k
&& \k $\smallmatrix \ovr{\text{Isomorphic}}\\ \text{fundamental}\\ \undr{\text{groups?}} \endsmallmatrix$ \k 
&& \k $\smallmatrix \text{$[L]$/$[L_c]$-isos-} \\ \text{pectral?} \endsmallmatrix$ \k
&& \k $\smallmatrix \text{$L$/$L_c$-isos-} \\ \text{pectral?} \endsmallmatrix$  \k
&& \k \text{Ex.} \quad \text{dim.} \k 
& \cr
\tablerule 
\k && $p$-iso., not $0$-iso.\  &&  No && Yes/No && No && No  
&& {\vrule height12pt width0pt depth8pt}  2.3 \quad $n\ge 4$  &\cr
&& all &&  Yes && No && No && Yes 
&& {\vrule height12pt width0pt depth8pt}  3.4 \quad $n \ge 4$  &\cr
&& $0$-iso., not $p$-iso.\ &&  No && No && No && Yes/No 
&& {\vrule height12pt width0pt depth8pt}  3.6  \quad $n\ge 6$  &\cr
 && all &&  Yes && No && Yes && Yes 
&& {\vrule height12pt width0pt depth8pt} 3.3 \quad $n\ge 4$  &\cr
&& all &&  Yes && Yes && Yes && Yes 
&& {\vrule height12pt width0pt depth8pt} 3.7 \quad $n\ge 7$  &\cr
&& none/some $p$-iso &&  No && No && Yes && Yes   
&& {\vrule height12pt width0pt depth8pt}  3.5 \quad $n\ge 13\!$ &\cr
 && none &&  No && No && No && Yes 
&& {\vrule height12pt width0pt depth8pt} 3.8 \quad $n\ge 4$  &\cr
\tablerule
\cr}}}

\mpb

We notice that all the pairs in the table above are not marked length isospectral.   
Some of the examples in the table 
have some similar spectral properties as other known examples 
in the context of nilmanifolds  (see \cite{Go} and \cite{Gt1}).  

The following proposition will show that marked length isospectral implies isometric for flat manifolds. 
This result adds more information to the table above.
The analogous result is known in other contexts, for instance, for flat tori, for closed surfaces of negative
curvature (see \cite{Ot} and \cite{Cr}) and for certain two-step nilmanifolds that include 
Heisenberg manifolds (see \cite{Eb}).

\proclaim{Proposition 3.10}
If two flat manifolds have the same marked length spectrum then they are isometric.
\endproclaim
\demo{Proof}
By assumption, there is an isomorphism, $\phi : \G \rightarrow\G'$, between the fundamental groups preserving
lengths, i.e.\ $l(\g)=l(\phi(\g))$, for any $\g\in\G$.
  
By Bieberbach's second theorem
$\phi$ is given by conjugation by $AL_a$, an affine motion. 
Since $\phi(L_\Ld)=AL_aL_\Ld{(AL_a)}^{-1}=L_{A\Ld}$, and on the other hand $\phi(L_\Ld)=L_{\Ld'}$ (since $\phi$ is an
isomorphism), it follows that 
$A\Ld=\Ld'$. Since $\phi$ preserves lengths, so does $A$,
hence $A\in O(n)$. Thus $\phi$ is given by conjugation by an isometry, hence $M_\G$ and $M_{\G'}$ are isometric.
\qed 
\enddemo

\heading{\S 4 Poisson summation formulas for flat manifolds}\endheading

We now consider, for $\G$ a Bieberbach group and $\tau$ a finite dimensional representation of $O(n)$, the zeta function
$$Z^{\Gamma}_\tau(s):= \sum_{\mu \ge0}\, d_{\tau,\mu}(\Gamma)\, e^{-4\pi^2 \mu s}.\tag{4.1}$$ 
The series is uniformly convergent for $s > \varepsilon$, for any $\varepsilon> 0$. 
We recall that for any $BL_b \in \G$ we have set $n_B=\text{dim}\, \text{ker}(B-\I)$ and $b_+ =p_B(b)$, 
where $p_B$ denotes the orthogonal projection onto $\text{ker}(B-\I)$.
In the case when $\tau = \tau_p$, for some 
$0\le p \le n$,  we write $Z^{\Gamma}_p(s)= Z^{\Gamma}_{\tau_p}(s)$.
\proclaim{Theorem 4.1} (i) We have
$$ Z^{\,\Gamma}_\tau(s)=
{|F|}^{-1} \sum_{BL_b \in F} \text{tr}\,\tau(B){\hbox{vol}({\Lambda ^*}^{\text{B}})}^{-1}\,
(4 \pi s)^{-\frac{n_B}2}
\sum _{\lambda_+ \in p_B(\Ld)} e^{- \frac{ {\| \lambda_+ + b_+\|}^2}{4s}}.\tag{4.2}
$$

(ii) $\text{Spec}(M_\Gamma)$ determines the lengths of closed geodesics of $M_\G$ and the numbers $n_B$.

(iii) $\text{Spec}_\tau(M_\Gamma)$ (in particular $\text{Spec}_p(M_\Gamma)$, for any $p\ge 0$)  determines 
the spectrum of the torus $T_\Ld=\Ld\backslash\R^n$ and the cardinality of $F$. That is, if $M_\G$ and $M_{\G'}$ are
$\tau$-isospectral then the associated tori $T_\Ld$ and $T_{\Ld'}$ are isospectral and $|F|=|F'|$.
\endproclaim

\demo{Proof}
Using expression (1.1) for $d_{\tau,\mu}(\Gamma)$  we may write
$$Z^{\Gamma}_\tau(s)= |F|^{-1} \sum_{BL_b \in F}\,\text{tr}\,\tau(B) 
\sum_{v \in {(\Lambda^*)}^{\text{B}}} e^{-2\pi i v.b}e^{-4\pi^2 \| v\|^2 s},\tag4.3$$ 
where ${(\Lambda^*)}^{\text{B}}:=   \Lambda^* \cap \text{ker}(B-\I)$. 
We have that ${(\Lambda^*)}^{\text{B}}$ is a lattice in $\text{ker}(B-\I)$, by 
Lemma 1.1(ii).

We shall recall some standard facts on Poisson summation.  
If $f \in {\Cal S}(\R^n)$, the Schwartz space of $\R^n$, let $\hat f (y) =\int_{\R^n} f(x) 
e^{-2\pi i x.y} dx$, the Fourier transform of $f$. We then have  (see \cite{Se}, for instance):
\roster
\item"--" If $h(x):= e^{-\pi \| x \|^2}$, then $\hat{h}=h$. 
If $a>0$ and  $g(x) := e^{-a\pi \| x \|^2}$ then
$\hat{g}(y)=a^{-\frac n2} h(\frac y{\sqrt a}).$
\item"--" If $b\in \R^n$ and $f \in {\Cal S}(\R^n)$, set  $f_b(x):= f(x) e^{2\pi
i x.b}$.   Then $\hat{f_b}(y)= \hat{f}(y-b)$.
\item"--"  If $L$ is a lattice in $\R^d$ and $L^*$ is the dual lattice of $L$ then: 
$$ \sum_{\nu\in L} f(\nu) = {\text{vol}(L)}^{-1}\sum_{\nu' \in L^*} \hat{f}(\nu').$$ 
\endroster

Now we apply Poisson summation in (4.3) for $Z^{\,\Gamma}_\tau(s)$, with $L={(\Lambda^*)}^{\text{B}}$,
a lattice in $\ker(B-\I)$, observing that in the expression for $e_{\mu,\g}$ in (1.1) we may write $b_+$ in place of
$b$.   By Lemma 1.1(iii) we get: 
$$ Z^{\,\Gamma}_\tau(s)=
{|F|}^{-1} \sum_{BL_b \in F} \text{tr}\,\tau(B){\text{vol} ({\Lambda ^*}^{\text{B}})}^{-1} 
\sum _{\lambda_+ \in p_B(\Ld)} (4 \pi s)^{-\frac{n_B}2} e^{- \frac{ \| \lambda_+ + b_+\|^2}{4s}}.
$$ 

Now, since $\text{tr}_0(B)=1$ for any $B$, the standard asymptotic argument implies that 
the lengths of closed geodesics $l(\gamma L_\lambda)=\| \lambda_+ + b_+\|$ and the numbers $n_B$, 
are determined by $\text{spec}(M_\Gamma)$ (see for instance \cite{Bu, \S 9.2}), hence (ii) follows. 

To verify (iii) we note first that
by (4.2) the $\tau$-spectrum of $M_\G$ determines the zeta function $Z^{\,\Gamma}_\tau(s)$. The standard asymptotic argument shows that this determines the following series

$$
{|F|}^{-1} \text{dim}(\tau)\, {\hbox{vol}(\Lambda ^*)}^{-1}\,
(4 \pi s)^{-\frac{n}2}
\sum _{\lambda \in \Ld} e^{- \frac{ {\| \lambda\|}^2}{4s}}$$
which is the partial sum of the right hand side of (4.2), corresponding  just to the element $BL_b=\I$ in $F$.
Now, by using Poisson summation for the torus, this expression is equal to
$${|F|}^{-1} \text{dim}(\tau)  \sum_{v\in {\Lambda^*}} e^{-4\pi^2 {\|v\|}^2 s}= {|F|}^{-1} \text{dim}(\tau)Z_1^\G(s).$$
Now we can leave out the factor $\text{dim}(\tau)$ and since the eigenvalue zero of the Laplacian on functions has
multiplicity one, $|F|$ is determined and hence the zeta function for the torus.
This completes the proof of the theorem.
\qed
\enddemo

\example{Remark 4.2} (a) The eigenvalue spectrum does not determine the complex 
lengths of closed geodesics, as Example 3.6 shows.

(b) Formula (4.2) is a Poisson summation formula for natural vector bundles over flat
manifolds. 
Sunada (see \cite{Su}) has obtained a similar formula in the case of functions, i.e.  $\tau$ is the trivial
representation, by using the heat kernel on $M_\Gamma$ and the Selberg trace formula. As a consequence, Sunada obtains
(iii) in the theorem in the function case.     The above approach is different since it  uses the formula for the
multiplicities of eigenvalues obtained in \cite{MR2} and furthermore the final formula is also different. 
\endexample

The $p^{\text{th}}$-Betti number of a closed Riemannian manifold gives the multiplicity of the eigenvalue zero for the Hodge-Laplacian acting on $p$-forms. 
Thus, a closed orientable $n$-manifold cannot be $n$-isospectral to a nonorientable one, since the $n^{th}$-Betti
numbers are distinct. In particular  such manifolds cannot be strongly isospectral.
One can ask whether
a closed orientable manifold can be isospectral on functions to a nonorientable one. As a
consequence of Theorem 4.1, we shall now show that this 
cannot happen for flat manifolds. We have:  

\noindent\proclaim{Corollary 4.3} 
If two flat manifolds are isospectral 
then they are both orientable or both nonorientable.
\endproclaim

\demo{Proof}
Let $M=M_\G$, $\G$ a Bieberbach group and let $\gamma=BL_b \in \Gamma$. The possible
eigenvalues of $B$ are 1,-1 and a set of complex roots of unity which come in pairs,
each one with the conjugate root. Hence, if $d_B^-$ denotes the multiplicity of the eigenvalue $-1$, 
$\text{det}(B)= (-1)^{d_B^-}= (-1)^{n- n_B}$.   Hence,  $M_\Gamma$ is orientable if and only if 
$n_B\equiv n$, mod $2$, for each $\gamma \in \Gamma$. On the other hand, all $n_B$'s are determined by 
the spectrum, by Theorem 4.1, thus the corollary follows.
\qed
\enddemo

\example{Remark 4.4}
The assertion in the corollary is not true for $p$-isospectral manifolds, $p>0$, as shown in Ex.\ 2.3(iii),(v) and in \cite{MR2,3}.
\endexample

If $\G$ is of diagonal type it is possible to
give a much more explicit formula for the zeta functions. 
Indeed, using the notations in (1.4) we have:

\proclaim{Theorem 4.5} 
If $\Gamma$ is a Bieberbach group of diagonal type with holonomy group isomorphic to $\Z_2^r$ then 
we have:
$$\text{(i)} \qquad\qquad\qquad Z^{\,\Gamma}_p(s)=
2^{-r} \sum _{d=1}^n K_p^n(n-d)\,
(4 \pi s)^{-\frac{d}2} \sum _{t=0}^d c_{d,t}(\G)\theta_{d,t}(\tfrac 1{4s})
\qquad \qquad\qquad\quad$$
where, for $\text{Re}\, z>0$,
$$ \theta_{d,t}(z) = \sum 
_{(m_1,\dots,m_d)\in\Z^d} 
 e^{-z \left(\sum_{j=1}^t (\frac 12 + m_j)^2 + \sum_{j=t+1}^d  {m_j}^2\right)}.$$ 

\noi (ii) If $\,\Gamma$ and $\Gamma'$ are Bieberbach groups of diagonal type with holonomy group $\Z_2^r$,
then $M_{\Gamma}$ and $M_{\Gamma'}$ are $p$-isospectral if and only if  
$$ K_p^n(n-d)\, c_{d,t}(\G)=K_p^n(n-d)\, c_{d,t}(\G')$$
for each $1\le t\le d \le n$.
In particular, if $c_{d,t}(\G)=c_{d,t}(\G')$ for every $d,t$, then
$M_\G$ and $M_{\G'}$ are $p$-isospectral for all $p$.

If $K_p^n(x)$ has no integral roots and $M_{\Gamma}$ and $M_{\Gamma'}$ 
are $p$-isospectral then they are Sunada isospectral. 
In particular, isospectrality implies Sunada isospectrality. 
\endproclaim

\demo{Proof}
In this case the volumes in (4.2) equal $1$. 
Using the notation $\tilde{I}_B= I_B\cap I_{2b}^{odd}$ and (1.4) 
we see that the zeta functions can be written:
$$
\align
 Z^{\,\Gamma}_p(s)&=
{|F|}^{-1} \sum_{BL_b \in F} K_p^n(n-n_B)\,
(4 \pi s)^{-\frac{n_B}2}
\!\!\sum _{m_j \in \Z : j \in I_B} \!\!
 e^{-\frac 1{4s} \left(\sum_{j\in \tilde{I}_B} (\frac 12 + m_j)^2 + \sum_{j \in I_B\setminus \tilde{I}_B} {m_j}^2\right)}\\
&=2^{-r} \sum _{d=1}^n K_p^n(n-d)\,
(4 \pi s)^{-\frac{d}2} \sum _{t=0}^d c_{d,t}(\G)
\!\!\!\!\sum_{(m_1,\dots,m_d)\in\Z^d}\!\!\!\!
 e^{-\frac 1{4s} \left(\sum_{j=1}^t (\frac 12 + m_j)^2 + \sum_{j=t+1}^d  {m_j}^2\right)}.
\endalign
$$
This implies the expression for $Z_p^{\,\G}(s)$ asserted in (i).

We furthermore note that, for each $d,t$, as $s\mapsto 0^+$,
$$(4 \pi s)^{-\frac{d}2} \theta_{d,t}(\tfrac 1{4s})\sim 2^t (4 \pi s)^{-\frac{d}2} e^{-\frac t{16s}}.
$$

Now,
$$ 
Z^{\,\Gamma}_p(s)\sim
2^{-r} \sum _{d=1}^n \sum _{t=0}^d K_p^n(n-d)\,
2^t  
c_{d,t}(\G)(4\pi s)^{-\frac{d}2} e^{-\frac t{16s}},\tag4.4
$$
as $s \mapsto 0^+$ and 
furthermore we have that $(4 \pi s)^{-\frac{d}2} e^{-\frac t{16s}}= o\left((4 \pi s)^{-\frac{d'}2} e^{-\frac
{t'}{16s}}\right)$ if and only if $t>t'$, or else $t=t'$ and $d<d'$. Thus, by a standard asymptotic argument we may
conclude from (4.4) that
if $Z^{\,\Gamma}_p(s)=Z^{\,\Gamma'}_p(s)$, then necessarily 
$K_p^n(n-d)\,c_{d,t}(\G)=K_p^n(n-d)\,c_{d,t}(\G')$, for every $d,t$. The converse is also clear from (4.4). 
This proves the first assertion in (ii).
Furthermore, if $K_p^n(x)$ has no integral roots, then  $p$-isospectrality implies  the equality of 
the numbers $c_{d,t}$ for $\G$ and $\G'$, hence  
$M_\G$ and $M_{\G'}$ are Sunada isospectral by Remark 1.3.
\qed
\enddemo

\example{Remark 4.6} 
The previous theorem says that, for groups of diagonal type, 
$p$-isos\-pec\-trality for one value of $p$ implies Sunada isospectrality, provided the Krawtchouk polynomial 
$K_p^n(x)$ has no integral roots. 
This extends to all dimensions a result obtained in \cite{MR3, Theorem 3.12(d)} for 
dimensions $\le 8$ ---using a very different approach.

We note that the theorem gives another proof of the fact that if $c_{d,t}(\G) = c_{d,t}(\G')$ for all $d,t$, then
 $M_\G$ and $M_{\G'}$ are $p$-isospectral for all $p$ (see Remark 1.3). 
\endexample

We furthermore have that the isospectrality criterion  in \cite{MR2, Thm.\ 3.1} is actually an equivalence, for
flat manifolds of diagonal type. That is, 

\proclaim{Proposition 4.7}
Let $\G$, $\G'$ be of diagonal type. 
Then $M_\G$ and $M_{\G'}$ are $p$-isospectral if and only if there is 
a bijection $\g\leftrightarrow \g'$ between the  holonomy groups $F$ and $F'$, 
such that for each $\mu$ $$\text{tr}_p(B) e_{\mu,\gamma}=\text{tr}_p(B') e_{\mu,\gamma'}$$ where $\g=BL_b$, $\g'=B'L_{b'}$.
In particular, $M_\G$ and $M_{\G'}$ are isospectral if and only if there is 
a bijection $\g\leftrightarrow \g'$ from $F$ onto $F'$ satisfying $e_{\mu,\gamma}=e_{\mu,\gamma'}$, 
for all $\g\in F$.
\endproclaim

\example{Remark}
We note that the bijection in the case of $p$-isospectrality is only necessary 
for those elements with nonvanishing trace. 
\endexample

\demo{Proof}
The ``if" part is stated in \cite{MR2, Thm.\ 3.1} and  follows directly from (1.1) with $\tau=\tau_p$. 
To prove the ``only if" part 
we set $F_d(\G)=\{\g=BL_b \in F: n_B=d \}$, defining $F'_d(\G')$ similarly for $\G'$. 
By Thm.\ 4.5 (ii),  $ K_p^n(n-d)\, c_{d,t}(\G)=K_p^n(n-d)\, c_{d,t}(\G')$ for every $d,t$, $t\le d$.
We distinguish two cases:

(a)  $d$ is such that $K_p^n(n-d)\ne 0, \quad$(b)  $d$ is such that $K_p^n(n-d)=0$.

For  $d$ verifying (a) we have that $c_{d,t}(\G)=c_{d,t}(\G')$ for every $t\le d$.
Hence, by (1.4) 
we can define a bijection $\g\leftrightarrow \g'$ between $F_d(\G)$ and $F_d(\G')$, such that $n_B=n_{B'}$ and  $|I_B\cap
I_{2b}^{odd}|=|I_{B'}\cap I_{2b'}^{odd}|$. 
Since $\G$ is of diagonal type, $e_{\mu,\gamma}$ depends only on $\mu$,
$n_B$ and $|I_B\cap I_{2b}^{odd}|$, 
hence this bijection is such that $e_{\mu,\gamma}=e_{\mu,\gamma'}$ for every $\g \in F_d(\G)$. 
Thus, there is bijection between the elements in $F$ and $F'$ with $n_B = n_B' =d$, $d$ verifying (a), such that  
$\text{tr}_p(B) e_{\mu,\gamma}=\text{tr}_p(B') e_{\mu,\gamma'}$. 

On the other hand, by Theorem 4.1,
$|F|=|F'|$, thus the partial bijection defined above implies that the number of elements $BL_b\in F$ with $n_B=d$ of type (b) (i.e. tr$_p(B)=0$) 
equals the number of elements $B'L_{b'}\in F'$
with $n_{B'}$ of type (b). For these elements, 
the equality $\text{tr}_p(B) e_{\mu,\gamma}=\text{tr}_p(B') e_{\mu,\gamma'}$ holds trivially, hence we can complete the bijection between $F$ and $F'$ as desired.
\qed
\enddemo

To conclude this paper, we will  compute explicitly the zeta functions $Z_p^\G(s)$ for some pairs of Bieberbach groups.  
We notice that from the expression of the zeta functions one can read $p$-isospectrality.
We consider first the case of the Klein bottle group. 

\example{Example 4.8} In the notation of Example 2.2 we have 
$$2Z_p^{\G}(s)= \tfrac {\binom 2p}{4 \pi s} \sum_{(m_1,m_2) \in \Z^2} e^{-\frac {{m_1}^2 + {m_2}^2 } {4s}} 
+ \tfrac {\text{tr}_p( B)}{\sqrt{4 \pi s}} \sum_{m\in \Z} e^{-\frac {\left(\frac 12 + m\right)^2 } {4s}}.\tag4.5
$$
We note furthermore that $\text{tr}_1( B)=0$, $\text{tr}_2( B)=-1$,  
thus from (4.5) we get an explicit formula for $Z_p^{\G}(s)$, $p=0,1,2$.
\endexample

\example{Example 4.9}
If we take $\G$ and $\G'$ as in Example 2.3(i),   
using the information therein, we get 
$$\align 2 \, Z_p^{\G}(s)&= \tfrac {\binom 4p}{(4 \pi s)^2} \sum_{(m_1,\dots,m_4) \in \Z^4} 
e^{-\frac {\sum_{j=1}^4 {m_j}^2} {4s}} 
+ \tfrac {K_p^4(3)}{\sqrt{4 \pi s}} \sum_{m\in \Z} e^{-\frac {\left(\frac 12 + m\right)^2 } {4s}}; \\
 2 \, Z_p^{\G'}(s)&= \tfrac {\binom 4p}{(4 \pi s)^2} \sum_{(m_1,\dots ,m_4) \in \Z^4} e^{-\frac {\sum_{j=1}^4 {m_j}^2} {4s}} 
+ \tfrac {K_p^4(1)}{(4 \pi s)^{\frac 32}} \sum_{(m_1,m_2,m_3) \in \Z^3} e^{-\frac {\left(\frac 12 + m_1\right)^2
+{m_2}^2+{m_3}^2} {4s}}.
\endalign$$
These expressions indicate that $\G$ and $\G'$ are $p$-isospectral if and only if $p=2$, 
since $K_2^4(1)=K_2^4(3)=0$ for $p=2$ only. 
Furthermore we also see that $M_\G$ and $M_{\G'}$ are neither isospectral nor $L$-isospectral.
\endexample

\example {Example 4.10} Using the information in Example 2.3(iii), in the case of $\G$ we have that
$n_{B_i}=2\,$ for $1\le i \le 3$ and in all three cases, $\text{tr}_p(B)= K_p^4(2)$ equals 0 if $p=1,3$, 
it equals $-2$ if $p=2$ and 1 if $p=0,4$. 

Thus, we obtain for $\G$:
$$ \alignat2
4Z_p^\G(s) &= \tfrac {\binom 4p}{(4 \pi s)^2} \sum_{(m_1,\dots ,m_4) \in \Z^4} 
e^{-\frac {\sum_{j=1}^4 {m_j}^2} {4s}} &&+ \tfrac {3\,K_p^4(2)}{4 \pi s} \sum_{(m_1,m_2) \in \Z^2} e^{-\frac
{\left(\frac 12 + m_1\right)^2 + {m_2}^2} {4s}}. 
\endalignat
$$ 

In the case of $\G'$, we have that
$n_{B'}=n_{{B'}^3}=1$, $n_{{B'}^2}=2$  and $\text{tr}_p(B')= \text{tr}_p({B'}^3)$ equals 0 if $p=1,2,3$ 
and equals 1 if $p=0,4$. Furthermore, $\text{tr}_p({B'}^2) = K_p^4(2)$. 
Thus we obtain:
$$
\alignat2 
4Z_p^{\G'}(s) &= \tfrac {\binom 4p}{(4 \pi s)^2} \sum_{(m_1,\dots ,m_4) \in \Z^4} 
e^{-\frac {\sum_{j=1}^4 {m_j}^2} {4s}}  &&+ \tfrac {K_p^4(2)}{4 \pi s} \sum_{(m_1,m_2) \in \Z^2} e^{-\frac
{\left(\frac 12 + m_1\right)^2 + {m_2}^2} {4s}}
\\ 
& &&+ \tfrac {\text{tr}_p(B')}{\sqrt{4 \pi s}} \sum_{m\in \Z}
\left(e^{-\frac {\left(\frac 14 + m\right)^2 } {4s}}+e^{-\frac {\left(\frac 34 + m\right)^2 } {4s}}\right).
\endalignat
$$
These expressions and the values of the $p$-traces 
imply that  $M_\G$ and $M_{\G'}$ are not $L$-isospectral and furthermore, that they are
$p$-isospectral if and only if  $p=1,3$.
Indeed, for these values of $p$ only, all $p$-traces are zero for any
$B\ne \I$ for both $\G$ and $\G'$, hence for such $p$ the expression of
$Z_p^{\,\G}(s)=Z_p^{\,\G'}(s)$ contains only the contribution of the lattice
elements.
\endexample

\example{Remark 4.11}
We now discuss how the $p$-isospectrality in Examples 2.3(i)-(vi), together with the the existence of different 
lengths of closed geodesics are  not in contradiction with  the wave trace and the heat trace formulas. The
basic point in both cases is that when considering the Laplacian acting on $p$-forms, certain coefficients 
appear multiplying the contributions of each geodesic to the formulas and they are not always positive, so they
can cancel out when added up, or they can vanish in some cases. The last situation happens in our examples.  

We first look at formula (4.2), that can be viewed as a heat 
trace formula
for the Laplacian on $p$-forms, concentrating on  Example 4.9 (i.e. Ex.\ 2.3(i)). In the comparison  between
$\G$ and $\G'$ we can see why the case $p=2$ is different. Indeed, in the expressions of the heat traces 
$Z_p^\Gamma(s)$ and $Z_p^{\Gamma'}(s)$, we find the coefficients $K_p^4(3)$ and $K_p^4(1)$ respectively  in the
right hand side, which vanish for $p=2$ and this makes both manifolds $2$-isospectral. The lengths of closed
geodesics are distinct for $M_\G$ and $M_{\G'}$ but there is no contradiction with the equality of the heat
traces for $p=2$, since these lengths show up in the exponents in the right hand side of the formula with 
coefficient  $0$ (for $p=2$), so they do not influence the sum. For other values of $p$ the coefficients
$K_p^4(3)$ and $K_p^4(1)$ do not vanish and we get quite different heat traces for $M_\G$ and $M_{\G'}$ showing
that the manifolds are not $p$-isospectral for $p\ne 2$.

A similar phenomenon happens with the singularities of the
wave traces, located at lengths of closed geodesics.

For the Laplacian acting on $p$-forms the following residue formula for the wave trace is stated in \cite{DG,
Introduction} under a genericity condition on the closed geodesics of period $T$: 
$$\lim_{t\to T} (t-T) \sum_{k\ge 0} e^{i \sqrt{\ld_k} \, t} =
\sum_{\g:l(\g)=T} \frac {T_{0_\g}}{2\pi} \, i^{\sigma_\g} \, {|I-P_\g|}^{-\frac 12}
\, \text{tr}(H_\g).\tag{4.6}$$
Here $T_{0_\g}$ is the smallest positive period of $\g$, $\sigma_\g$ is 
the Maslov factor, $P_\g$  the Poincar\'e map around $\g$ and $H_\g:\Ld^p\mapsto\Ld^p$ 
the holonomy along $\g$. We observe in (4.6) that the factor tr$(H_\g)$ is in our flat case the $p$-trace of
the orthogonal transformation $B$ denoted by tr$_p(B)$    
which vanishes for some values of $p$ (as we have just seen) 
and,  in this case,  this implies the vanishing of the r.h.s.\ in (4.6).
For instance in Example 2.3(i), $\frac{\sqrt{5}}2$ is the length of a closed geodesic
for $M_{\G'}$ but not for $M_\G$. For $M_{\G'}$ it is a singularity of the wave trace 
in (4.6) for the Laplacian acting on functions, but it is not a singularity for the wave trace of
the Laplacian acting on 2-forms since tr$_2(B')=0$, 
for $\g'=B'L_{b'} \in \G'$ corresponding to the geodesics of length $\frac{\sqrt{5}}2$ in $M_{\G'}$.
\endexample

\enddocument